\documentclass{article}
\usepackage{amsmath}
\usepackage{amssymb}
\usepackage{amsthm}
\usepackage{graphicx}
\usepackage{epic,eepic,epsfig}
\usepackage{color}
\usepackage{caption}
\usepackage{subcaption}
\usepackage{placeins}	% add 1.22
\usepackage{multirow}
\usepackage{epstopdf}
\usepackage{varwidth}
\usepackage[table]{xcolor}

\newtheorem{theorem}{Theorem}[section]
\newtheorem{lemma}[theorem]{Lemma}

\theoremstyle{definition}

\newcommand{\rev}[1]{#1}

\newcommand{\bsym}[1]{\mathbf{#1}}

\definecolor{tabclr}{cmyk}{0,0,1,0}

\begin{document}

\title{An embedded SDG method for the convection-diffusion equation}

\author{
Siu Wun Cheung\thanks{Department of Mathematics, Texas A\&M University, College Station, TX 77843, USA (\texttt{tonycsw2905@math.tamu.edu})}
\and
Eric T. Chung\thanks{Department of Mathematics, The Chinese University of Hong Kong, Shatin, New Territories, Hong Kong SAR, China (\texttt{tschung@math.cuhk.edu.hk}) }
}

\maketitle

\begin{abstract}
In this paper, we present an embedded staggered discontinuous Galerkin method for the convection-diffusion equation. The new method combines the advantages of staggered discontinuous Galerkin (SDG) and embedded discontinuous Galerkin (EDG) method, and results in many good properties, namely local and global conservations, free of carefully designed stabilization terms or flux conditions and high computational efficiency. In applying the new method to convection-dominated problems, the method provides optimal convergence in potential and suboptimal convergence in flux, which is comparable to other existing DG methods, and achieves $L^2$ stability by making use of a skew-symmetric discretization of the convection term, irrespective of diffusivity. We will present numerical results to show the performance of the method.
\end{abstract}

\section{Introduction}
\label{sec:intro}

Discontinuous Galerkin (DG) methods were first introduced by Reed and Hill for solving hyperbolic equations \cite{first-dg}. The DG methods have been proven superior to the classical continuous Galerkin (CG) methods for hyperbolic problems. In the past two decades, DG methods have also been applied to second-order elliptic problems. A comprehensive study on DG methods for elliptic problems is given in \cite{unified-dg}. The original DG methods for elliptic problems, using polynomial approximations of degree $k$ for both the potential and the flux, converge with optimal order $k+1$ for the potential but suboptimal order $k$ for the flux. While the same orders of convergence can be obtained by using classical CG finite element methods, the DG methods give rise to a discrete problem with a higher number of degrees of freedom. DG methods have therefore been criticized for its high computation cost and judged to be not being particularly useful for elliptic problems. Later, the hybridizible discontinuous Galerkin (HDG) method was introduced for solving elliptic problems \cite{unified-hdg}. The HDG method provides optimal orders of convergence for both the potential and the flux in $L^2$ norm. Moreover, superconvergence can be obtained for the potential through a local postprocessing technique. 

In recent years, there are active developments of DG methods for problems in fluid dynamics and wave propagations, see for example \cite{carrero05,hdg,cockburn05,ldg,conv-diff,houston09,Liu-Shu,hpdg,shahbazi07,nguyen11}. 
On the other hand, staggered meshes bring the advantages of reducing numerical dissipation in computational fluid dynamics \cite{sfvm,sdm,sdm1},
and numerical dispersion in computational wave propagation \cite{semi, fully, newdg, newdg1, meta, jcp-max,geo}.
Combining the ideas of DG methods and staggered meshes, a new class of staggered discontinuous Galerkin (SDG) methods was proposed for approximations of Stokes system \cite{kcl-2013-Stokes-SDG}, convection-diffusion equation \cite{convdiff}, and incompressible Navier-Stokes equations \cite{sdg-ns1}. The new class of SDG methods possesses many good properties, including local and global conservations, stability in energy, and optimal convergence. For a more complete discussion on the SDG method, see also \cite{newdg,newdg1,meta,jcp-max,curlcurl,convdiff,kcl-2013-FETI-DP-Stokes} and the references therein.

In \cite{sdg-hdg, sdg-hdg1}, it was shown that the SDG method can be regarded as a limit of the HDG method. The SDG method can be obtained from the HDG method by setting the stabilization parameter on a set of edges to be zero and letting the parameter on another set of edges to infinity. As a result, the SDG method inherits the advantages of the HDG method, including superconvergence through the use of a local postprocessing technique. Furthermore, in the SDG method for incompressible Navier-Stokes equations \cite{sdg-ns1}, using the postprocessing and a spectro-consistent discretizations with a novel splitting of the diffusion and the convection term, stability in $L^2$ energy is achieved. 

The embedded discontinuous Galerkin (EDG) method was first introduced for solving the linear shell problems \cite{edg-shell}. Later, an EDG method for solving second order elliptic problems was discussed and analyzed in \cite{edg-elliptic}. The EDG method was obtained from HDG method by enforcing strong continuity for hybrid unknowns \cite{unified-hdg}. This greatly reduces the number of degrees of freedom in the globally coupled system and makes the EDG method has a higher computational efficiency compared with other DG methods. As a tradeoff for this advantage, the EDG method is not locally conservative and loses the optimal convergence in the flux achieved by the HDG method \cite{edg-cfd}. The loss in accuracy makes the EDG method a less attractive candidate compared with the HDG method. However, the optimal order of convergence for HDG method is also lost in the case of convection-dominated problems as shown the numerical examples \cite{hdg-convdiff}. In this case, the EDG method becomes appealing alternative to all other DG methods including the HDG method, since it has a higher computational efficiency and the same orders of convergence. On the other hand, compared with the CG finite element method, the EDG method provides the same sparsity structure of the stiffness matrix after static condensation, whilist the EDG method is more robust, accurate and stable than the CG finite element method in convection-dominated problems. We remark that the multiscale discontinuous Galerkin (MDG) method \cite{buffa06, hughes06} are related to the EDG method, which is originally proposed for the convection-diffusion problems. The MDG method and the EDG method are both designed for a globally continuous approximation of the solution and can give rise to identical schemes. Recently, the EDG method has been proposed on Euler equations and Navier-Stokes equations \cite{edg-cfd, edg-ns}. Due to the advantages shared with other DG methods and the high computational efficiency compared with other DG methods, the EDG method has also been applied to challenging problems in computational fluid dynamics, such as implicit large eddy simulation \cite{edg-les}. 

In this paper, we propose a combination of the SDG method and the EDG method for the convection-diffusion equation. The new method seeks approximations in the SDG locally conforming finite element spaces, which gives rise to a flux formulation without introducing any carefully designed stabilization terms or flux conditions as in other DG methods. The new method further reduces the size of the global discrete problem compared with SDG method by restricting the numerical approximation for the primal unknown to lie in a proper subspace of the SDG finite element space. Moreover, the new method inherits the stability in $L^2$ energy thanks to spectro-consistent discretizations in the SDG method. The convergence is optimal with order $k+1$ for \rev{the unknown function} and suboptimal with order $k$ for \rev{its gradient}, which are comparable to all other DG methods for convection-dominated problems. 

The paper is organized as follows. In Section \ref{sec:esdg}, we will have a derivation on the method. Next, in Section \ref{sec:stab}, we will provide a stability analysis of the method. Then, in Section \ref{sec:num}, we will present extensive numerical examples to see the performance of our method. Finally, a conclusion is given.

\section{Method description}
\label{sec:esdg}

\subsection{Model problem}
\label{sec:prob}
Let $\Omega \subset \mathbb{R}^2$ be a polygonal domain. We consider the steady-state convection-diffusion equation with homogeneous Dirichlet boundary condition:
\begin{equation}
\begin{split}
- \mu \bsym{\Delta} u + \text{div} \, (\bsym{b}u) & = f \; \mbox{ in }\Omega,\\
u & = 0 \; \mbox { on } \partial \Omega.\\
\end{split}
\label{eq:convdiff1}
\end{equation}
Here $u$ is the unknown function to be approximated, $\bsym{b} = (b_1, b_2)$ is a divergence-free convection field, $f$ is a given source term and $g$ is a given boundary condition. Also, $\mu$ is the diffusivity, which is assumed to be constants throughout the domain $\Omega$. 
Before we start the derivation of our method, we shall state the variational formulation of the problem. Suppose $\bsym{b} \in L^{\infty}(\Omega)$ and $f \in H^{-1}(\Omega)$. The variational formulation of the convection-diffusion equation is given by: find $u \in H_0^1(\Omega)$ such that for any $v \in H_0^1(\Omega)$, we have
\begin{equation}
\mu (\nabla u, \nabla v)_{0,\Omega} + (\text{div} \, (\bsym{b}u), v)_{0,\Omega} = (f, v)_{0,\Omega}.
\label{eq:convdiff_weak}
\end{equation}
Here $(\cdot, \cdot)_{0,\Omega}$ denotes the standard $L^2(\Omega)$ inner product. 

We will derive a mixed method for the problem. 
Since $\text{div} \, \bsym{b} = 0$, it is direct to see that
\begin{equation}
\begin{split}
\text{div} \, (\bsym{b}u) = \bsym{b} \cdot \nabla u + (\text{div} \, \bsym{b}) u = \bsym{b} \cdot \nabla u.
\end{split}
\label{eq:adjoint_cont}
\end{equation}
We can therefore rewrite \eqref{eq:convdiff1} as
\begin{equation}
\begin{split}
- \mu \bsym{\Delta} u + \dfrac{1}{2} \text{div} \, (\bsym{b}u) + \dfrac{1}{2} \bsym{b} \cdot \nabla u & = f \; \mbox{ in }\Omega,\\
u & = 0 \; \mbox { on } \partial \Omega.\\
\end{split}
\label{eq:convdiff_split}
\end{equation}
We introduce the auxiliary variables
\begin{equation}
\begin{split}
\bsym{w} & = \sqrt{\mu}  \, \nabla u - \frac{1}{2\sqrt{\mu}} \bsym{b} u,\\
\bsym{p} & = \bsym{b} u.\\
\end{split}
\label{eq:newvar}
\end{equation}
Then \eqref{eq:convdiff_split} can be reformulated as a system of first-order linear PDEs:
\begin{equation}
\begin{split}
-\sqrt{\mu} \, \text{div}\, \bsym{w} + \frac{1}{2\sqrt{\mu}} \bsym{b}\cdot\bsym{w} + \frac{1}{4 \mu} \bsym{b}\cdot\bsym{p} & = f \; \mbox{ in }\Omega, \\
u & = 0 \; \mbox{ on }\partial \Omega.
\end{split}
\label{eq:convdiff2}
\end{equation}

\subsection{Staggered meshes}
\label{sec:mesh}

\begin{figure}[ht!]
\centering
\includegraphics[width=3.0in]{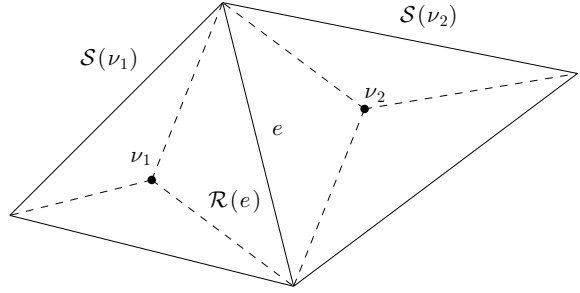}
\caption{An illustration of the staggered mesh in two dimensions.}\label{fig_mesh}
\end{figure}

Let $\mathcal{T}_u$ be a triangulation of the two-dimensional domain $\Omega$ by a
set of triangles without hanging nodes.
We introduce the notation $\mathcal{F}_u$ to denote the
set of all edges in the triangulation $\mathcal{T}_u$ and
$\mathcal{F}_u^0$ to denote the subset of all interior edges in
$\mathcal{F}_u$ excluding those on the boundary of $\Omega$.
We also denote the set of all vertices in $\mathcal{T}_u$ by $\mathcal{N}_u$.
For each triangle in $\mathcal{T}_u$, we take an interior point $\nu$, denote the initial triangle by $\mathcal{S}(\nu)$, and divide $\mathcal{S}(\nu)$ into three triangles by
joining the point $\nu$ and the three vertices of $\mathcal{S}(\nu)$.
We also denote the set of all interior points $\nu$ by $\mathcal{N}$, the set of all new edges generated by the subdivision of triangles by $\mathcal{F}_p$, and the triangulation after subdivision by $\mathcal{T}$.
Note that the interior point $\nu$ of each triangle in $\mathcal{T}_u$ should be chosen such that the new triangulation $\mathcal{T}$ observes the shape regularity criterion.
In practice, we can simply choose $\nu$ as the centroid of the triangle.
Also,
$\mathcal{F} = \mathcal{F}_u\cup\mathcal{F}_p$
denotes the set of all edges of triangles in $\mathcal{T}$
and $\mathcal{F}^0 = \mathcal{F}^0_u \cup \mathcal{F}_p$ denotes the set of all interior edges of triangles in $\mathcal{T}$.
For each edge $e \in \mathcal{F}_u$, we let $\mathcal{R}(e)$ be the
union of the all triangles in the new triangulation $\mathcal{T}$ sharing the edge $e$.
{Figure}~\ref{fig_mesh} demonstrates these definitions. The edges $e \in \mathcal{F}_u$ are represented in solid lines and the $e \in \mathcal{F}_p$ are represented in dotted lines.
%Finally, for $\eta \in \mathcal{N}_u$, we let $\mathcal{Q}(\eta)$ be the union of all initial triangles $\mathcal{S}(\nu)$ with a vertex $\eta$.

For each edge $e \in \mathcal{F}$, we will also define a unit normal vector $\bsym{n}_{e}$ in the following way.
If $e \in \mathcal{F} \setminus \mathcal{F}^0$ is a boundary edge, then we define $\bsym{n}_{e}$ as the outward unit normal vector of $e$ from $\Omega$.
If $e\in\mathcal{F}^0$ is an interior edge, then $\bsym{n}_{e}$ is fixed as one of the two possible unit normal vectors on $e$.
When it is clear that which edge we are considering, we omit the index $e$ and write the unit normal vector as $\bsym{n}$.

To end this section, we define the jumps in the following way:
for any edge $e \in \mathcal{F}$, denote one of the triangles in the refined triangulation $\mathcal{T}$, which contains $e$ by $\tau^+$, and denote the other triangle, if exists, by $\tau^-$.
The outward unit normal vectors on $e$ in $\tau^+$ and $\tau^-$ are denoted by $\bsym{n}^+$ and $\bsym{n}^-$, respectively.
Also, for any quantity $\phi$, the notations $\phi^\pm$ are defined on the edge $e$ by the values of $\phi \vert_{\tau^{\pm}}$ restricted on $e$.
Then, if $\phi$ is a scalar quantity, the notation $[\phi]$ over an edge $e$ defined as
\begin{equation}
[\phi]\vert_e:= (\bsym{n} \cdot \bsym{n}^+) \phi^+ + (\bsym{n} \cdot \bsym{n}^-) \phi^-.
\end{equation}
If $\bsym{\Phi}$ is a vector quantity, then the notation $[\bsym{\Phi} \cdot \bsym{n}]$ is similarly defined as
\begin{equation}
[\bsym{\Phi} \cdot \bsym{n}]|_e:=(\bsym{n}\cdot \bsym{n}^+) (\bsym{\Phi}^+ \cdot \bsym{n})
+ (\bsym{n}\cdot \bsym{n}^-) (\bsym{\Phi}^- \cdot \bsym{n}).
\end{equation}

\subsection{SDG and ESDG finite element spaces}
\label{sec:SDGspace}

We will define the finite element spaces.
Let $k \geq 0$ be a non-negative integer.
Let $\tau\in\mathcal{T}$ and $e\in\mathcal{F}$.
We define $P^{k}(\tau)$ and $P^k(e)$ as the space of polynomials whose order is not greater than $k$ on $\tau$ and $e$, respectively.
We will also define norms on the spaces.
We use the standard notations $ \| \cdot \|_{0, \Omega}$ to denote the standard $L^2$ norm on $\Omega$ and $ \| \cdot \|_{0, e}$ to denote the $L^2$ norm on an edge $e$.

First, we define the following locally $H^{1}(\Omega)$-conforming finite element space:
\begin{equation}
U^{h}
= \{ v \: : \: v|_{\tau} \in P^{k}(\tau); \; \tau\in\mathcal{T}; \;
v \;\text{is
continuous over} \; e \in \mathcal{F}_{u}^{0}; \;
v|_{\partial \Omega}=0 \}.
\label{eq:Uh}
\end{equation}
Note that for any $v\in U^h$, we have $v|_{\mathcal{R}(e)} \in
H^1(\mathcal{R}(e))$ for each edge $e \in \mathcal{F}_u$.
We define the following discrete $L^2$-norm $\| \cdot \|_X$ and discrete $H^1$-norm $\| \cdot \|_Z$ on the space $U^h$:
\begin{equation}
\begin{split}
\| v \|_X & = \left( \| v \|_{0, \Omega}^2 + \sum_{e \in \mathcal{F}_u^0} h_e \| v \|_{0,e}^2 \right)^\frac{1}{2}, \\
\| v \|_Z & = \left( \| \nabla_h v \|_{0, \Omega}^2 + \sum_{e \in \mathcal{F}_p} h_e^{-1} \| [v] \|_{0,e}^2 \right)^\frac{1}{2}.
\end{split}
\label{eq:Unorm}
\end{equation}

Next, we define the following ESDG finite element space, which is a proper subspace of the SDG finite element space:
\begin{equation}
\widetilde{U}^{h}
= \{ v \in U^h \: : \: 
v \;\text{is
continuous at} \; \eta; \; \eta \in \mathcal{N}_u \}.
\label{eq:tildeUh}
\end{equation}
Note that the test functions $v \in \widetilde{U}^h$ are continuous at only all the nodes in the initial grid but not the nodes at the refined grid, are therefore they are not globally continuous. We remark that for the EDG method \cite{edg-elliptic}, the space of the numerical trace is imposed with a global continuity on the skeleton of the mesh. 

Finally, we define the following locally $H(\text{div};\Omega)$-conforming finite element space:
\begin{equation}
W^{h}
= \{ \bsym{\Psi} \: : \: \bsym{\Psi}|_{\tau} \in P^{k}(\tau)^{2}; \: \tau\in\mathcal{T}; \:
\bsym{\Psi} \cdot \bsym{n} \: \text{is continuous over} \:
e\in \mathcal{F}_{p} \}.
\label{eq:Wh}
\end{equation}
Note that for any $\bsym{\Psi} \in W^h$, we have
$\bsym{\Psi} \vert_{\mathcal{S}(\nu)} \in H(\text{div};\mathcal{S}(\nu))$
for each $\nu \in \mathcal{N}$.
We define the following discrete $L^2$-norm $\| \cdot \|_{X'}$ and discrete $H(\text{div};\Omega)$-norm $\| \cdot \|_{Z'}$ on the space $W^h$:
\begin{equation}
\begin{split}
\| \bsym{\Psi} \|_{X'} & = \left( \| \bsym{\Psi} \|_{0, \Omega}^2 + \sum_{e \in \mathcal{F}_p} h_e \| \bsym{\Psi} \cdot \bsym{n} \|_{0,e}^2 \right)^\frac{1}{2}, \\
\| \bsym{\Psi} \|_{Z'} & = \left( \| \text{div}_h \bsym{\Psi} \|_{0, \Omega}^2 + \sum_{e \in \mathcal{F}_u^0} h_e^{-1} \| [\bsym{\Psi} \cdot \bsym{n}] \|_{0,e}^2 \right)^\frac{1}{2}.
\end{split}
\label{eq:Wnorm}
\end{equation}
\rev{By \cite{newdg,newdg1}, there exist interpolation operators $\mathcal{I}$ onto $\widetilde{U}^h$ and $\mathcal{J}$ onto $W^h$ such that
\begin{equation}
\begin{split}
B_h^*(u - \mathcal{I} u,  \bsym{\Psi}) & = 0 \text{ for all } \bsym{\Psi} \in W^h, \\
B_h(\mathbf{w} - \mathcal{J} \mathbf{w}, v) & = 0 \text{ for all } v \in U^h,
\label{eq:ortho_esdg}
\end{split}
\end{equation}
and
\begin{equation}
\begin{split}
\| u - \mathcal{I} u \|_{0,\Omega} & \leq Ch^{k+1} \vert u \vert_{H^{k+1}(\Omega)}, \\
\| \mathbf{w} - \mathcal{J} \mathbf{w} \|_{0,\Omega} & \leq Ch^{k+1} \vert \mathbf{w} \vert_{[H^{k+1}(\Omega)]^2}.
\end{split}
\label{eq:approx_esdg}
\end{equation}
}

\subsection{Derivation of the method}
\label{sec:sdg}

We will derive the discrete problem in our SDG formulation
starting from the system of first order equations in \eqref{eq:newvar} and \eqref{eq:convdiff2}.

Multiplying the first equation of \eqref{eq:newvar} by $\bsym{\Psi}_1 \in W^h$
and integrating over $\mathcal{S}(\nu)$ for $\nu \in \mathcal{N}$, we obtain
\begin{equation}
\int_{\mathcal{S}(\nu)} \bsym{w} \cdot \bsym{\Psi}_1 \; dx = - \sqrt{\mu} \int_{\mathcal{S}(\nu)} u (\text{div } \bsym{\Psi}_1) \; dx
+ \sqrt{\mu} \int_{\partial \mathcal{S}(\nu)} u (\bsym{\Psi}_1 \cdot \bsym{n}) \; d\sigma
- \frac{1}{2\sqrt{\mu}} \int_{\mathcal{S}(\nu)} \bsym{p} \cdot \bsym{\Psi}_1 \; dx.
\label{eq:w1}
\end{equation}
Similarly, multiplying the second equation of \eqref{eq:newvar} by $\bsym{\Psi}_2 \in W^h$
and integrating over $\mathcal{S}(\nu)$ for $\nu \in \mathcal{N}$, we have
\begin{equation}
\int_{\mathcal{S}(\nu)} \bsym{p} \cdot \bsym{\Psi}_2 \; dx = \int_{\mathcal{S}(\nu)} u (\bsym{b} \cdot \bsym{\Psi}_2) \; dx.
\label{eq:tw1}
\end{equation}

Finally, multiplying the first equation of \eqref{eq:convdiff2} by $v \in U^h$
and integrating over $\mathcal{R}(e)$ for $e \in \mathcal{F}_u^0$, we have
\begin{equation}
\begin{split}
& \sqrt{\mu} \int_{\mathcal{R}(e)} \bsym{w} \cdot \nabla v\; dx - \sqrt{\mu} \int_{\partial \mathcal{R}(e)} ( \bsym{w}\cdot \bsym{n}) v \;d\sigma
+ \frac{1}{2\sqrt{\mu}} \int_{\mathcal{R}(e)} (\bsym{b}\cdot \bsym{w}) \, v \; dx \\
& + \frac{1}{4 \mu} \int_{\mathcal{R}(e)} (\bsym{b}\cdot \bsym{p}) \, v \; dx
= \int_{\mathcal{R}(e)} f v\; dx.
\end{split}
\label{eq:div_w}
\end{equation}

Summing those equations in \eqref{eq:w1}--\eqref{eq:div_w} over all $\mathcal{R}(e)$ and $\mathcal{S}(\nu)$,
we obtain the staggered discontinuous Galerkin method for
\eqref{eq:convdiff1} proposed in \cite{convdiff}: find
$(u_h,\bsym{w}_h,\bsym{p}_h)\in U^h \times W^h \times W^h$ such that
for any $v \in U^h, \bsym{\Psi}_1, \bsym{\Psi}_2 \in W^h$, we have
\begin{equation}
\begin{split}
\sqrt{\mu} B_h(\bsym{w}_h,v)+\frac{1}{2\sqrt{\mu}}R_h\left(\bsym{w}_h + \frac{1}{2\sqrt{\mu}}\bsym{p}_h,v\right) & = (f,v)_{0,\Omega},\\
\sqrt{\mu} B_h^*(u_h,\bsym{\Psi}_1) - \frac{1}{2\sqrt{\mu}} (\bsym{p}_h,\bsym{\Psi}_1)_{0,\Omega} & = (\bsym{w}_h,\bsym{\Psi}_1)_{0,\Omega},\\
R_h^*(u_h,\bsym{\Psi}_2) &=(\bsym{p}_h,\bsym{\Psi}_2)_{0,\Omega},\\
\end{split}
\label{eq:convdiff_sdg}
\end{equation}
where bilinear forms $B_h(\bsym{\Psi},v)$ and $B^*_h(v,\bsym{\Psi})$ are defined as
\begin{equation}
\begin{split}
B_{h}(\bsym{\Psi},v) &=  \int_{\Omega} \bsym{\Psi} \cdot \nabla_h v
\; dx
- \sum_{e\in\mathcal{F}_{p}}\int_{e} (\bsym{\Psi} \cdot \bsym{n}) \: [v] \; d\sigma,  \\
B^{*}_{h}(v,\bsym{\Psi}) &= -\int_{\Omega}  v \: \text{div}_h\,
\bsym{\Psi} \;dx + \sum_{e\in\mathcal{F}_{u}^{0}}\int_{e} v \:
[\bsym{\Psi} \cdot \bsym{n}] \; d\sigma,
\end{split}
\label{eq:bilinear1}
\end{equation}
and the bilinear forms $R_h(\bsym{\Psi},v)$ and $R_h^*(v,\bsym{\Psi})$ are defined as
\begin{equation}
\begin{split}
R_h(\bsym{\Psi},v) & = \int_\Omega (\bsym{b} \cdot \bsym{\Psi}) \, v \; dx,\\
R_h^*(v,\bsym{\Psi}) & = \int_\Omega v \, (\bsym{b} \cdot \bsym{\Psi}) \; dx.
\end{split}
\label{eq:bilinear3}
\end{equation}

Now, if we consider only the test functions $v \in \widetilde{U}^h$ in \eqref{eq:div_w}, and we seek an approximation $\widetilde{u}_h \in \widetilde{U}^h$ for the unknown function $u$, we obtain the embedded staggered discontinuous Galerkin method for
\eqref{eq:convdiff1}: find
$(\widetilde{u}_h,\widetilde{\bsym{w}}_h,\widetilde{\bsym{p}}_h)\in \widetilde{U}^h \times W^h \times W^h$ such that
for any $v \in \widetilde{U}^h, \bsym{\Psi}_1, \bsym{\Psi}_2 \in W^h$, we have
\begin{equation}
\begin{split}
\sqrt{\mu} B_h(\widetilde{\bsym{w}}_h,v)+\frac{1}{2\sqrt{\mu}}R_h\left(\widetilde{\bsym{w}}_h + \frac{1}{2\sqrt{\mu}}\widetilde{\bsym{p}}_h,v\right) & = (f,v)_{0,\Omega},\\
\sqrt{\mu} B_h^*(\widetilde{u}_h,\bsym{\Psi}_1) - \frac{1}{2\sqrt{\mu}} (\widetilde{\bsym{p}}_h,\bsym{\Psi}_1)_{0,\Omega} & = (\widetilde{\bsym{w}}_h,\bsym{\Psi}_1)_{0,\Omega},\\
R_h^*(\widetilde{u}_h,\bsym{\Psi}_2) &=(\widetilde{\bsym{p}}_h,\bsym{\Psi}_2)_{0,\Omega}.\\
\end{split}
\label{eq:convdiff_esdg}
\end{equation}

By \cite{newdg1}, the two bilinear forms in \eqref{eq:bilinear1} satisfy the adjoint relation
\begin{equation}
B_h(\bsym{\Psi},v) = B_h^*(v,\bsym{\Psi})
\label{eq:adjoint-B}
\end{equation}
for all $v\in U_h$ and $\bsym{\Psi} \in W^h$.
The bilinear forms $B_h$ and $B_h^*$ are also continuous with respect to suitable discrete norms
\begin{equation}
\begin{split}
\vert B_h(\bsym{\Psi}, v) \vert & \leq \| \bsym{\Psi} \|_{X'} \| v \|_Z, \\
\vert B_h^*(v, \bsym{\Psi}) \vert & \leq  \| v \|_X \| \bsym{\Psi} \|_{Z'}, \\
\end{split}
\label{eq:cont-B}
\end{equation}
for all $v\in U_h$ and $\bsym{\Psi} \in W^h$.
Moreover, the bilinear forms $B_h$ and $B_h^*$ satisfy a pair of inf-sup conditions: there exists constants $\beta_1$ and $\beta_2$, independent of $h$, such that
\begin{equation}
\begin{split}
\inf_{v \in U^h \setminus \{0\}} \sup_{\bsym{\Psi} \in W^h \setminus \{\bsym{0}\}}
\frac{B_h(\bsym{\Psi}, v)}{ \| \bsym{\Psi} \|_{X'} \| v \|_Z} & \geq \beta_1, \\
\inf_{\bsym{\Psi} \in W^h \setminus \{\bsym{0}\}} \sup_{v \in U^h \setminus \{0\}}
\frac{B_h^*(v, \bsym{\Psi})}{\| v \|_X \| \bsym{\Psi} \|_{Z'}} & \geq \beta_2. \\
\end{split}
\label{eq:inf-sup-B}
\end{equation}

Also, it is obvious that the two bilinear forms in (\ref{eq:bilinear3}) satisfy
\begin{equation}\label{eq:adjoint-r}
R_h^*(v,\bsym{\Psi})=R_h(\bsym{\Psi},v)
\end{equation}
for all $v \in U_h$ and $\bsym{\Psi} \in W^h$.

\subsection{Linear system}
\label{sec:sys}

In this section,
we derive the linear systems resulting from \eqref{eq:convdiff_sdg} and \eqref{eq:convdiff_esdg}. We denote the corresponding matrix representation of the bilinear forms $B_h$ and $R_h$ by $B$ and $R$, respectively. Then by the adjoint properties, the matrix representation of the bilinear forms $B^*_h$ and $R^*_h$ are given by $B^T$ and $R^T$, respectively. Also, the notations for the finite element solutions would be abused to denote their corresponding vector representations.

Using these notations, we can write the SDG method \eqref{eq:convdiff_sdg} as a linear system of algebraic equations. The second equation of \eqref{eq:convdiff_sdg} can be written as
\begin{equation}
\sqrt{\mu} B^T u_h - \frac{1}{2\sqrt{\mu}} M\bsym{p}_h = M\bsym{w}_h,
\label{eq:matrix_eq1}
\end{equation}
where $M$ is the mass matrix for the space $W^h$.
Similarly, the last equation of \eqref{eq:convdiff_sdg} can be written as
\begin{equation}
R^T u_{h,1} = M\bsym{p}_h.
\label{eq:matrix_eq2}
\end{equation}
Lastly, the first equations of \eqref{eq:convdiff_sdg} can be written as
\begin{equation}
\sqrt{\mu} B\bsym{w}_h + \frac{1}{2\sqrt{\mu}} R\left(\bsym{w}_h + \frac{1}{2\sqrt{\mu}}\bsym{p}_h\right)
= f_h.
\label{eq:matrix_eq3}
\end{equation}
We can now obtain a linear system with the unknowns $\bsym{w}_h$ and $\bsym{p}_h$ eliminated. Combining \eqref{eq:matrix_eq1} and \eqref{eq:matrix_eq2}, we have
\begin{equation}
\begin{split}
\bsym{w}_h & = M^{-1} \left(\sqrt{\mu} B^T u_{h} - \frac{1}{2\sqrt{\mu}} R^T u_{h}\right),\\
\bsym{p}_h & = M^{-1} R^T u_{h}.
\end{split}
\label{eq:matrix_eq4}
\end{equation}
We note that the elimination can be done by solving small problems in each $\mathcal{S}(\nu)$ since $M$ is a block diagonal matrix with each block corresponding to the mass matrix
of $W^h|_{\mathcal{S}(\nu)}$.

We further introduce the notations
\begin{equation}
\begin{split}
\bsym{\Delta}_h &= -BM^{-1}B^T,\\
\bsym{b} \cdot \nabla_h &= -\frac{1}{2} BM^{-1}R^T + \frac{1}{2} RM^{-1}B^T, \\
A &= -\mu \bsym{\Delta}_h + \bsym{b} \cdot \nabla_h.
\end{split}
\label{eq:mat_def}
\end{equation}
We note that the discrete diffusion operator $-\Delta_h$ is symmetric and positive-definite, and the discrete convection operator $\bsym{b} \cdot \nabla_h$ is skew-symmetric.
Combining \eqref{eq:matrix_eq3} and \eqref{eq:matrix_eq4}, the algebraic system of the discrete problem \eqref{eq:convdiff_sdg} can then be reduced to
\begin{equation}
A u_h = f_h.
\label{eq:convdiff_sdg_system}
\end{equation}

Now, if we denote the matrix representation of the canonical embedding $\iota: \widetilde{U}^h \to U^h$ by $P$, then the matrix representations $\widetilde{B}$ and $\widetilde{R}$ of the bilinear forms $B_h \vert _{W^h \times \widetilde{U}^h}$ and $R_h \vert _{W^h \times \widetilde{U}^h}$ are related to $B$ and $R$ by
\begin{equation}
\begin{split}
\widetilde{B} &= PB,\\
\widetilde{R} &= PR. \\
\end{split}
\label{eq:mat_rep_esdg}
\end{equation}
The corresponding matrix represenations for the discrete diffusion operator and discrete convection operator are
\begin{equation}
\begin{split}
\widetilde{\bsym{\Delta}}_h &= -\widetilde{B}M^{-1}\widetilde{B}^T = P \bsym{\Delta}_h P^T,\\
\bsym{b} \cdot \widetilde{\nabla}_h &= -\frac{1}{2} \widetilde{B}M^{-1}\widetilde{R}^T + \frac{1}{2} \widetilde{R}M^{-1}\widetilde{B}^T = P(\bsym{b} \cdot \nabla_h) P^T. \\
\end{split}
\label{eq:mat_def_esdg}
\end{equation}
Therefore the algebraic system of the discrete problem \eqref{eq:convdiff_esdg} is given by
\begin{equation}
\widetilde{A} \widetilde{u}_h = \widetilde{f}_h,
\label{eq:convdiff_esdg_system}
\end{equation}
where
\begin{equation}
\begin{split}
\widetilde{A} & = P A P^T, \\
\widetilde{f}_h & = P f_h.
\end{split}
\end{equation}
We remark that in the embedded SDG method, the discretization of the diffusion operator $-\widetilde{\bsym{\Delta}}_h$ is still symmetric and positive definite. Similarly, the discretization of the convection operator $\bsym{b} \cdot \widetilde{\nabla}_h$ is still skew-symmetric. Therefore the spectro-consistent discretization is preserved in the new method.

\section{Stability analysis}
\label{sec:stab}

We start this section by stating a stability result in $L^2$ energy for the variational formulation \eqref{eq:convdiff_weak}.
\begin{lemma}
\label{lemma:stab_weak}
Let $u \in H_0^1(\Omega)$ be the weak solution of the variational formulation \eqref{eq:convdiff_weak} of the convection-diffusion equation.
Then we have 
\begin{equation}
\mu \| \nabla u \|_{0,\Omega}^2 = (f, u)_{0,\Omega}.
\end{equation}
\begin{proof}
In \eqref{eq:convdiff_weak}, we take a test function $v = u$. Then we have
\begin{equation}
\mu \| \nabla u \|_{0,\Omega}^2 + (\text{div} \, (\bsym{b} u), u)_{0,\Omega} = (f, u)_{0,\Omega}.
\end{equation}
On the other hand, using an integration by parts and the relation \eqref{eq:adjoint_cont}, we have
\begin{equation}
\begin{split}
(\text{div} \, (\bsym{b} u), u)_{0,\Omega}
& = - (\bsym{b}u, \nabla u)_{0,\Omega} \\
& = - (\bsym{b} \cdot \nabla u, u)_{0,\Omega} \\
& = - (\text{div} \, (\bsym{b} u), u)_{0,\Omega},
\end{split}
\end{equation}
which implies $(\text{div} \, (\bsym{b} u), u)_{0,\Omega} = 0$. This completes the proof.
\end{proof}
\end{lemma}

We will next see that the ESDG method provides a similar stability result. In the SDG method, the stability in $L^2$ energy is due to a spectro-consistent discretizations with the splitting of the diffusion and the convection term proposed in \cite{convdiff}. The stability in $L^2$ energy in a numerical method for the convection-diffusion problems is a kind of measure of how well the numerical solution approximates the analytical solution, and has significant effects on the quality of the numerical solution (see, for example, \cite{sdg-ns1}, \cite{sdg-ibm}; also see Section \ref{sec:exp3}). The ESDG method inherits the stability in $L^2$ energy from the SDG method due to the same spectro-consistent discretization structure.

The unknowns in the space $W^h$ in both the SDG method and the ESDG method give rise to an approximation of the flux $\bsym{z} = \nabla u$ in the space $W^h$. 
For the ESDG method, suppose $(\widetilde{u}_h,\widetilde{\bsym{w}}_h,\widetilde{\bsym{p}}_h)\in \widetilde{U}^h \times W^h \times W^h$ is the solution of \eqref{eq:convdiff_esdg}. An approximation $\widetilde{\bsym{z}}_h \in W^h$ for the flux $\bsym{z}$ is then given by
\begin{equation}
\widetilde{\bsym{z}}_h = \frac{1}{\sqrt{\mu}} \widetilde{\bsym{w}}_h + \frac{1}{2 \mu} \widetilde{\bsym{p}}_h = M^{-1} B^T P^T \widetilde{u}_h.
\label{eq:flux_esdg}
\end{equation}
Likewise for the SDG method, suppose $(u_h,\bsym{w}_h,\bsym{p}_h)\in U^h \times W^h \times W^h$ is the solution of \eqref{eq:convdiff_sdg}. An approximation $\bsym{z}_h \in W^h$ for the flux $\bsym{z}$ is given by
\begin{equation}
\bsym{z}_h = \frac{1}{\sqrt{\mu}} \bsym{w}_h + \frac{1}{2 \mu} \bsym{p}_h = M^{-1} B^T u_h.
\label{eq:flux_sdg}
\end{equation}

We are now ready to state the stability result for the ESDG method:
\begin{lemma}
\label{lemma:stab_esdg}
Let $(\widetilde{u}_h,\widetilde{\bsym{w}}_h,\widetilde{\bsym{p}}_h)\in \widetilde{U}^h \times W^h \times W^h$ be the numerical solution of the ESDG method \eqref{eq:convdiff_esdg}.
Then we have 
\begin{equation}
\mu \| \widetilde{\bsym{z}}_h \|_{0,\Omega}^2 = (f, \widetilde{u}_h)_{0,\Omega},
\end{equation}
where $\widetilde{\bsym{z}}_h \in W^h$ is defined in \eqref{eq:flux_esdg}.
\begin{proof}
In \eqref{eq:convdiff_esdg}, we take test functions as follows:
\begin{equation}
\begin{split}
v & = \widetilde{u}_h, \\
\bsym{\Psi}_1 & = -\widetilde{\bsym{w}}_h \\
\bsym{\Psi}_2 & = -\dfrac{1}{2} \widetilde{\bsym{z}}_h.
\end{split}
\end{equation}
Then we have
\begin{equation}
\begin{split}
\sqrt{\mu} B_h(\widetilde{\bsym{w}}_h,\widetilde{u}_h)+\frac{1}{2\sqrt{\mu}}R_h\left(\widetilde{\bsym{w}}_h + \frac{1}{2\sqrt{\mu}}\widetilde{\bsym{p}}_h,\widetilde{u}_h\right) & = (f,\widetilde{u}_h)_{0,\Omega},\\
-\sqrt{\mu} B_h^*(\widetilde{u}_h,\widetilde{\bsym{w}}_h) + \frac{1}{2\sqrt{\mu}} (\widetilde{\bsym{p}}_h,\widetilde{\bsym{w}}_h)_{0,\Omega} & = -(\widetilde{\bsym{w}}_h,\widetilde{\bsym{w}}_h)_{0,\Omega},\\
-\dfrac{1}{2}R_h^*(\widetilde{u}_h,\widetilde{\bsym{z}}_h) & = - \dfrac{1}{2}(\widetilde{\bsym{p}}_h,\widetilde{\bsym{z}}_h)_{0,\Omega}.\\
\end{split}
\label{eq:stab3.2.1}
\end{equation}
Recalling the definition of $\widetilde{\bsym{z}}_h$ in \eqref{eq:flux_esdg}, the above equations can be rewritten as
\begin{equation}
\begin{split}
\sqrt{\mu} B_h(\widetilde{\bsym{w}}_h,\widetilde{u}_h)+\frac{1}{2}R_h\left(\widetilde{\bsym{z}}_h,\widetilde{u}_h\right) & = (f,\widetilde{u}_h)_{0,\Omega},\\
-\sqrt{\mu} B_h^*(\widetilde{u}_h,\widetilde{\bsym{w}}_h) + \sqrt{\mu} (\widetilde{\bsym{z}}_h,\widetilde{\bsym{w}}_h)_{0,\Omega} & = 0,\\
-\dfrac{1}{2}R_h^*(\widetilde{u}_h,\widetilde{\bsym{z}}_h) + \dfrac{1}{2}(\widetilde{\bsym{z}}_h,\widetilde{\bsym{p}}_h)_{0,\Omega} & = 0.\\
\end{split}
\label{eq:stab3.2.2}
\end{equation}
Summing up the equations in \eqref{eq:stab3.2.2}, using the adjoint relations \eqref{eq:adjoint-B} and \eqref{eq:adjoint-r}, and the definition of $\widetilde{\bsym{z}}_h$ in \eqref{eq:flux_esdg} again, we have
\begin{equation}
\begin{split}
\sqrt{\mu} (\widetilde{\bsym{z}}_h,\widetilde{\bsym{w}}_h)_{0,\Omega} + \dfrac{1}{2}(\widetilde{\bsym{z}}_h,\widetilde{\bsym{p}}_h)_{0,\Omega} & = (f,\widetilde{u}_h)_{0,\Omega} \\
\mu \left(\widetilde{\bsym{z}}_h, \dfrac{1}{\sqrt{\mu}}\widetilde{\bsym{w}}_h + \dfrac{1}{2 \mu} \widetilde{\bsym{p}}_h\right)_{0,\Omega} & = (f,\widetilde{u}_h)_{0,\Omega} \\
\mu \| \widetilde{\bsym{z}}_h \|_{0,\Omega}^2 & = (f,\widetilde{u}_h)_{0,\Omega}.
\end{split}
\label{eq:stab3.2.3}
\end{equation}
\end{proof}
\end{lemma}
An important message from Lemma \ref{lemma:stab_esdg} is that the convection field $\bsym{b}$ vanishes in the above $L^2$ stability estimate for $\widetilde{\bsym{z}}_h$. This makes the ESDG approximation mimics the weak solution $\bsym{z}$ better as the convection field $\bsym{b}$ also vanishes in the above $L^2$ stability estimate for $\bsym{z}$ in Lemma \ref{lemma:stab_weak}. This is an advantage brought by the novel splitting of the convection term and the diffusion term.

To end this section, we establish the main stability result.
\begin{theorem}
\label{thm:stab_maon}
Let $(\widetilde{u}_h,\widetilde{\bsym{w}}_h,\widetilde{\bsym{p}}_h)\in \widetilde{U}^h \times W^h \times W^h$ be the numerical solution of the ESDG method \eqref{eq:convdiff_esdg}.
Then we have 
\begin{equation}
\mu \| \widetilde{u}_h \|_{Z} \leq C \| f \|_{0,\Omega},
\end{equation}
where $C$ is a constant independent of mesh size and diffusivity.
\begin{proof}
By Lemma \ref{lemma:stab_esdg}, we have
\begin{equation}
\mu \| \widetilde{z}_h \|_{0,\Omega}^2 = (f, \widetilde{u}_h)_{0,\Omega}.
\end{equation}
By Cauchy-Schwarz inequality and the equivalence of the standard $L^2$ norm $\| \cdot \|_{0,\Omega}$ and the discrete $H^1$ norm $\| \cdot \|_Z$ on the finite dimensional space $\widetilde{U}^h$, we obtain the following estimate for the right hand side:
\begin{equation}
(f, \widetilde{u}_h)_{0,\Omega} \leq \| f \|_{0,\Omega} \|\widetilde{u}_h\|_{0,\Omega} \leq K \| f \|_{0,\Omega} \|\widetilde{u}_h\|_Z,
\end{equation}
where $K$ is the constant from the equivalence of norms.
On the other hand, by the adjoint relation \eqref{eq:adjoint-B} and the first inf-sup condition in \eqref{eq:inf-sup-B}, we have
\begin{equation}
\begin{split}
\| \widetilde{u}_h \|_Z  
& \leq \dfrac{1}{\beta_1} \sup_{\bsym{\Psi} \in W^h \setminus \{\bsym{0}\}} \dfrac{B_h^*(\widetilde{u}_h, \bsym{\Psi})}{\| \bsym{\Psi} \|_{X'}}\\
& \leq \dfrac{1}{\beta_1} \sup_{\bsym{\Psi} \in W^h \setminus \{\bsym{0}\}} \dfrac{B_h^*(\widetilde{u}_h, \bsym{\Psi})}{\| \bsym{\Psi} \|_{0,\Omega}} \\
& = \dfrac{1}{\beta_1} \sup_{\bsym{\Psi} \in W^h \setminus \{\bsym{0}\}} \dfrac{(\widetilde{\bsym{z}}_h, \bsym{\Psi})_{0,\Omega}}{\| \bsym{\Psi} \|_{0,\Omega}} \\
& = \dfrac{1}{\beta_1} \| \widetilde{\bsym{z}}_h \|_{0,\Omega}.
\end{split}
\label{eq:energy-inf-sup}
\end{equation}
Therefore we have
\begin{equation}
\mu \| \widetilde{u}_h \|_Z^2 \leq \beta_1^2 (f, \widetilde{u}_h)_{0,\Omega} \leq \beta_1^2 K \| f \|_{0,\Omega} \|\widetilde{u}_h\|_Z.
\end{equation}
Dividing $\| \widetilde{u}_h \|_Z$ on both sides, we obtain the desired result.
\end{proof}
\end{theorem}

\section{Convergence analysis}
\label{sec:conv}
\rev{In this section, we present an error estimate between the weak solution $u$ in \eqref{eq:convdiff2} and the ESDG solution $\widetilde{u}_h$ in \eqref{eq:convdiff_esdg}.
\begin{theorem}
Let $(u,\bsym{w},\bsym{p})$ be the solution of \eqref{eq:newvar}--\eqref{eq:convdiff2}. Let $(\widetilde{u}_h,\widetilde{\bsym{w}}_h,\widetilde{\bsym{p}}_h)\in \widetilde{U}^h \times W^h \times W^h$ be the numerical solution of the ESDG method \eqref{eq:convdiff_esdg}.
Then we have the following optimal error bound:
\begin{equation}
\| u - \widetilde{u}_h \|_{L^2(\Omega)} \leq C(1+\mu^{-1})h^{k+1},
\end{equation}
where $C$ is a constant independent of mesh size and diffusivity.
\begin{proof}
First, we note that the solution $(u,\bsym{w},\bsym{p})$ satisfies the following system:
\begin{equation}
\begin{split}
\sqrt{\mu} B_h(\bsym{w},v)+\frac{1}{2\sqrt{\mu}}R_h\left(\bsym{w} + \frac{1}{2\sqrt{\mu}}\bsym{p},v\right) & = (f,v)_{0,\Omega},\\
\sqrt{\mu} B_h^*(u,\bsym{\Psi}_1) - \frac{1}{2\sqrt{\mu}} (\bsym{p},\bsym{\Psi}_1)_{0,\Omega} & = (\bsym{w}_h,\bsym{\Psi}_1)_{0,\Omega},\\
R_h^*(u,\bsym{\Psi}_2) &=(\bsym{p},\bsym{\Psi}_2)_{0,\Omega}.\\
\end{split}
\label{eq:convdiff_reform}
\end{equation}
for any $v \in \widetilde{U}^h, \bsym{\Psi}_1, \bsym{\Psi}_2 \in W^h$.
Subtracting \eqref{eq:convdiff_esdg} from \eqref{eq:convdiff_reform}, we have
\begin{equation}
\begin{split}
\sqrt{\mu} B_h(\bsym{w}-\widetilde{\bsym{w}}_h,v)+\frac{1}{2\sqrt{\mu}}R_h\left((\bsym{w} - \widetilde{\bsym{w}}_h) + \frac{1}{2\sqrt{\mu}}(\bsym{p}-\widetilde{\bsym{p}}_h),v\right) & = 0,\\
\sqrt{\mu} B_h^*(u-\widetilde{u}_h,\bsym{\Psi}_1) - \frac{1}{2\sqrt{\mu}} (\bsym{p}-\widetilde{\bsym{p}}_h,\bsym{\Psi}_1)_{0,\Omega} & = (\bsym{w}-\widetilde{\bsym{w}}_h,\bsym{\Psi}_1)_{0,\Omega} \\
R_h^*(u - u_h,\bsym{\Psi}_2) &=(\bsym{p} - \widetilde{\bsym{p}}_h,\bsym{\Psi}_2)_{0,\Omega}.\\
\end{split}
\label{eq:convdiff_differ1}
\end{equation}
Introduce the notations
\begin{equation}
\begin{split}
\delta_u = \mathcal{I}u - \widetilde{u}_h \in \widetilde{U}^h, \quad &
\quad \varepsilon_u = u - \mathcal{I}u, \\
\bsym{\delta}_w = \mathcal{J}\bsym{w} - \widetilde{\bsym{w}}_h \in W^h, \quad & 
\quad \bsym{\varepsilon}_w = \bsym{w} - \mathcal{J}\bsym{w}, \\
\bsym{\delta}_p = \mathcal{J}\bsym{p} - \widetilde{\bsym{p}}_h \in W^h, \quad &
\quad \bsym{\varepsilon}_p = \bsym{p} - \mathcal{J}\bsym{p}.\\
\end{split}
\end{equation}
Using the properties in \eqref{eq:ortho_esdg}, we can rewrite \eqref{eq:convdiff_differ1} as
\begin{equation}
\begin{split}
\sqrt{\mu} B_h(\bsym{\delta}_w,v)+\frac{1}{2\sqrt{\mu}}R_h\left(\bsym{\delta}_w + \frac{1}{2\sqrt{\mu}}\bsym{\delta}_p,v\right) & = -\frac{1}{2\sqrt{\mu}}R_h\left(\bsym{\varepsilon}_w + \frac{1}{2\sqrt{\mu}}\bsym{\varepsilon}_p,v\right),\\
\sqrt{\mu} B_h^*(\delta_u,\bsym{\Psi}_1) - \left(\bsym{\delta}_w + \frac{1}{2\sqrt{\mu}} \bsym{\delta}_p,\bsym{\Psi}_1\right)_{0,\Omega} & = \left(\bsym{\varepsilon}_w + \frac{1}{2\sqrt{\mu}} \bsym{\varepsilon}_p,\bsym{\Psi}_1\right)_{0,\Omega},\\
R_h^*(\delta_u,\bsym{\Psi}_2) -(\bsym{\delta}_p,\bsym{\Psi}_2)_{0,\Omega} & = -R_h^*(\varepsilon_u,\bsym{\Psi}_2) +(\bsym{\varepsilon}_p,\bsym{\Psi}_2)_{0,\Omega}.\\
\end{split}
\label{eq:convdiff_differ2}
\end{equation}
Using the argument as in \eqref{eq:energy-inf-sup}, by the second equation in \eqref{eq:convdiff_differ2}, we have
\begin{equation}
\| \delta_u \|_Z  \leq \dfrac{1}{\beta_1 \sqrt{\mu}} \left(\left\| \bsym{\delta}_w + \frac{1}{2\sqrt{\mu}} \bsym{\delta}_p \right\|_{0,\Omega} + \left\| \bsym{\varepsilon}_w + \frac{1}{2\sqrt{\mu}} \bsym{\varepsilon}_p \right\| \right).
\end{equation}
Moreover, using a discrete Poincar\'{e} inequality, we have
\begin{equation}
\| \delta_u \|_{0,\Omega}  \leq \dfrac{K}{\beta_1 \sqrt{\mu}} \left(\left\| \bsym{\delta}_w + \frac{1}{2\sqrt{\mu}} \bsym{\delta}_p \right\|_{0,\Omega} + \left\| \bsym{\varepsilon}_w + \frac{1}{2\sqrt{\mu}} \bsym{\varepsilon}_p \right\|_{0,\Omega} \right).
\label{eq:conv4.1.0}
\end{equation}
On the other hand, in \eqref{eq:convdiff_differ2}, we take we take test functions as follows:
\begin{equation}
\begin{split}
v & = \delta_u, \\
\bsym{\Psi}_1 & = -\bsym{\delta}_w, \\
\bsym{\Psi}_2 & = -\dfrac{1}{2\sqrt{\mu}} \bsym{\delta}_w - \dfrac{1}{4\mu} \bsym{\delta}_p.
\end{split}
\end{equation}
Then we have
\begin{equation}
\begin{split}
\sqrt{\mu} B_h(\bsym{\delta}_w,\delta_u)+\frac{1}{2\sqrt{\mu}}R_h\left(\bsym{\delta}_w + \frac{1}{2\sqrt{\mu}}\bsym{\delta}_p,\delta_u\right) & = -\frac{1}{2\sqrt{\mu}}R_h\left(\bsym{\varepsilon}_w + \frac{1}{2\sqrt{\mu}}\bsym{\varepsilon}_p,\delta_u\right),\\
-\sqrt{\mu} B_h^*(\delta_u,\bsym{\delta}_w) + \left(\bsym{\delta}_w + \frac{1}{2\sqrt{\mu}} \bsym{\delta}_p,\bsym{\delta}_w\right)_{0,\Omega} & = -\left(\bsym{\varepsilon}_w + \frac{1}{2\sqrt{\mu}} \bsym{\varepsilon}_p,\bsym{\delta}_w\right)_{0,\Omega},\\
 -\dfrac{1}{2\sqrt{\mu}}R_h^*\left(\delta_u,\bsym{\delta}_w+\dfrac{1}{2\sqrt{\mu}}\bsym{\delta}_p\right) + \dfrac{1}{2\sqrt{\mu}}\left(\bsym{\delta}_p,\bsym{\delta}_w+\dfrac{1}{2\sqrt{\mu}}\bsym{\delta}_p\right)_{0,\Omega} & =  \dfrac{1}{2\sqrt{\mu}}R_h^*\left(\varepsilon_u,\bsym{\delta}_w+\dfrac{1}{2\sqrt{\mu}}\bsym{\delta}_p\right) \\ 
 & \quad - \dfrac{1}{2\sqrt{\mu}}\left(\bsym{\varepsilon}_p,\bsym{\delta}_w+\dfrac{1}{2\sqrt{\mu}}\bsym{\delta}_p\right)_{0,\Omega}.\\
\end{split}
\label{eq:conv4.1.1}
\end{equation}
Summing up the equations in \eqref{eq:conv4.1.1} and using the adjoint relations \eqref{eq:adjoint-B} and \eqref{eq:adjoint-r}, we have
\begin{equation}
\left\| \bsym{\delta}_w+\dfrac{1}{2\sqrt{\mu}}\bsym{\delta}_p \right\|_{0,\Omega}^2 = T_1 + T_2 + T_3 + T_4,
\label{eq:conv4.1.2}
\end{equation}
where
\begin{equation}
\begin{split}
T_1 & = \dfrac{1}{2\sqrt{\mu}}R_h^*\left(\varepsilon_u,\bsym{\delta}_w+\dfrac{1}{2\sqrt{\mu}}\bsym{\delta}_p\right), \\
T_2 & = - \dfrac{1}{2\sqrt{\mu}}\left(\bsym{\varepsilon}_p,\bsym{\delta}_w+\dfrac{1}{2\sqrt{\mu}}\bsym{\delta}_p\right)_{0,\Omega}, \\
T_3 & = -\frac{1}{2\sqrt{\mu}}R_h\left(\bsym{\varepsilon}_w + \frac{1}{2\sqrt{\mu}}\bsym{\varepsilon}_p,\delta_u\right), \\
T_4 & = -\left(\bsym{\varepsilon}_w + \frac{1}{2\sqrt{\mu}} \bsym{\varepsilon}_p,\bsym{\delta}_w\right)_{0,\Omega}.
\end{split}
\end{equation}
Next, we will estimate each of these terms. Using Young's inequality, we have
\begin{equation}
\begin{split}
\vert T_1 \vert & \leq \dfrac{1}{2\sqrt{\mu}} \| \bsym{b} \|_{L^\infty(\Omega)} \|\varepsilon_u\|_{0,\Omega} \left\|\bsym{\delta}_w + \frac{1}{2\sqrt{\mu}}\bsym{\delta}_p\right\|_{0,\Omega} \\
& \leq \dfrac{1}{4\mu} \| \bsym{b} \|_{L^\infty(\Omega)}^2 \|\varepsilon_u\|_{0,\Omega}^2 + \dfrac{1}{4} \left\|\bsym{\delta}_w + \frac{1}{2\sqrt{\mu}}\bsym{\delta}_p\right\|_{0,\Omega}^2.
\end{split}
\label{eq:conv4.1.3a}
\end{equation}
Similarly, for $T_2$, we imply
\begin{equation}
\begin{split}
\vert T_2 \vert & \leq  \dfrac{1}{2\sqrt{\mu}} \|\bsym{\varepsilon}_p\|_{0,\Omega} \left\|\bsym{\delta}_w + \frac{1}{2\sqrt{\mu}}\bsym{\delta}_p\right\|_{0,\Omega} \\
& \leq \dfrac{1}{4\mu} \|\bsym{\varepsilon}_p\|_{0,\Omega}^2 + \dfrac{1}{4}  \left\|\bsym{\delta}_w + \frac{1}{2\sqrt{\mu}}\bsym{\delta}_p\right\|_{0,\Omega}^2.
\end{split}
\label{eq:conv4.1.3b}
\end{equation}
For $T_3$, we have
\begin{equation}
\begin{split}
\vert T_3 \vert & \leq \dfrac{1}{2\sqrt{\mu}} \| \bsym{b} \|_{L^\infty(\Omega)} \left\|\bsym{\varepsilon}_w + \frac{1}{2\sqrt{\mu}}\bsym{\varepsilon}_p\right\|_{0,\Omega} \|\delta_u\|_{0,\Omega}, \\
& \leq \dfrac{2K^2}{\beta_1^2 \mu^2} \| \bsym{b} \|_{L^\infty(\Omega)}^2 \left\|\bsym{\varepsilon}_w + \frac{1}{2\sqrt{\mu}}\bsym{\varepsilon}_p\right\|_{0,\Omega}^2 + \dfrac{\beta_1^2 \mu}{32 K^2}\|\delta_u\|_{0,\Omega}^2. \\
\end{split}
\label{eq:conv4.1.3c}
\end{equation}
For $T_4$, we first observe that
\begin{equation}
\begin{split}
\vert T_4 \vert & \leq \left\| \bsym{\varepsilon}_w + \frac{1}{2\sqrt{\mu}} \bsym{\varepsilon}_p \right\|_{0,\Omega} \left \| \bsym{\delta}_w \right\|_{0,\Omega} \\
& \leq \left\| \bsym{\varepsilon}_w + \frac{1}{2\sqrt{\mu}} \bsym{\varepsilon}_p \right\|_{0,\Omega} \left(\left\| \bsym{\delta}_w + \dfrac{1}{2\sqrt{\mu}} \bsym{\delta}_p \right\|_{0,\Omega} + \dfrac{1}{2\sqrt{\mu}} \left\| \bsym{\delta}_p \right\|_{0,\Omega} \right).
\end{split}
\end{equation}
Taking $\bsym{\Psi}_2 = -\bsym{\delta}_p$ in the last equation of \eqref{eq:convdiff_differ2}, we have
\begin{equation}
\begin{split}
\| \bsym{\delta}_p \|_{0,\Omega}^2 
& = R_h^*(\delta_u, \bsym{\delta}_p) + R_h^*(\varepsilon_u, \bsym{\delta}_p) - (\bsym{\varepsilon}_p, \bsym{\delta}_p)_{0,\Omega} \\
& \leq \| \bsym{b} \|_{L^\infty(\Omega)} (\| \delta_u \|_{0,\Omega} + \| \varepsilon_u \|_{0,\Omega}) \| \bsym{\delta}_p \|_{0,\Omega} + \| \bsym{\varepsilon}_p \|_{0,\Omega} \| \bsym{\delta}_p \|_{0,\Omega}
\end{split}
\end{equation}
Hence we imply
\begin{equation}
\begin{split}
\vert T_4 \vert & \leq \left\| \bsym{\varepsilon}_w + \frac{1}{2\sqrt{\mu}} \bsym{\varepsilon}_p \right\|_{0,\Omega} \left \| \bsym{\delta}_w \right\|_{0,\Omega} \\
& \leq \left\| \bsym{\varepsilon}_w + \frac{1}{2\sqrt{\mu}} \bsym{\varepsilon}_p \right\|_{0,\Omega} \left(\left\| \bsym{\delta}_w + \dfrac{1}{2\sqrt{\mu}} \bsym{\delta}_p \right\|_{0,\Omega} + \dfrac{1}{2\sqrt{\mu}} \left(  \| \bsym{b} \|_{L^\infty(\Omega)} (\| \delta_u \|_{0,\Omega} + \| \varepsilon_u \|_{0,\Omega}) + \| \bsym{\varepsilon}_p \|_{0,\Omega} \right) \right) \\
& \leq \left\| \bsym{\varepsilon}_w + \frac{1}{2\sqrt{\mu}} \bsym{\varepsilon}_p \right\|_{0,\Omega}^2 + \dfrac{1}{4}\left\| \bsym{\delta}_w + \dfrac{1}{2\sqrt{\mu}} \bsym{\delta}_p \right\|_{0,\Omega}^2 + \dfrac{2K^2}{\beta_1^2 \mu^2} \| \bsym{b} \|_{L^\infty(\Omega)}^2 \left\|\bsym{\varepsilon}_w + \frac{1}{2\sqrt{\mu}}\bsym{\varepsilon}_p\right\|_{0,\Omega}^2 + \dfrac{\beta_1^2 \mu}{32 K^2}\|\delta_u\|_{0,\Omega}^2 \\
& \quad \quad + \left\| \bsym{\varepsilon}_w + \frac{1}{2\sqrt{\mu}} \bsym{\varepsilon}_p \right\|_{0,\Omega}^2 + \dfrac{1}{16\mu} \left(  \| \bsym{b} \|_{L^\infty(\Omega)}^2 \| \varepsilon_u \|_{0,\Omega}^2 + \| \bsym{\varepsilon}_p \|_{0,\Omega}^2 \right) 
\end{split}
\label{eq:conv4.1.3d}
\end{equation}
Combining all these estimates with \eqref{eq:conv4.1.2}, we have
\begin{equation}
\dfrac{1}{4} \left\| \bsym{\delta}_w+\dfrac{1}{2\sqrt{\mu}}\bsym{\delta}_p \right\|_{0,\Omega}^2
\leq C \left(\mu^{-1} \| \varepsilon_u \|_{0,\Omega}^2 + \mu^{-1} \| \bsym{\varepsilon}_p \|_{0,\Omega}^2 + (1+\mu^{-2}) \left\| \bsym{\varepsilon}_w + \dfrac{1}{2\sqrt{\mu}} \bsym{\varepsilon}_p \right\|_{0,\Omega}^2 \right)  + \dfrac{\beta_1^2 \mu}{16 K^2} \|\delta_u\|_{0,\Omega}^2.
\label{eq:conv4.1.4}
\end{equation}
Combining \eqref{eq:conv4.1.0} and \eqref{eq:conv4.1.4}, we have
\begin{equation}
\begin{split}
\| \delta_u \|_{0,\Omega}^2 
& \leq \dfrac{K^2}{\beta_1^2 \mu} \left(\left\| \bsym{\delta}_w + \frac{1}{2\sqrt{\mu}} \bsym{\delta}_p \right\|_{0,\Omega} + \left\| \bsym{\varepsilon}_w + \frac{1}{2\sqrt{\mu}} \bsym{\varepsilon}_p \right\|_{0,\Omega} \right)^2 \\
&  \leq \dfrac{2 K^2}{\beta_1^2 \mu} \left(\left\| \bsym{\delta}_w + \frac{1}{2\sqrt{\mu}} \bsym{\delta}_p \right\|_{0,\Omega}^2 + \left\| \bsym{\varepsilon}_w + \frac{1}{2\sqrt{\mu}} \bsym{\varepsilon}_p \right\|_{0,\Omega}^2 \right) \\
& \leq \dfrac{2 K^2}{\beta_1^2 \mu} \left(C \left(\mu^{-1} \| \varepsilon_u \|_{0,\Omega}^2 + \mu^{-1} \| \bsym{\varepsilon}_p \|_{0,\Omega}^2 + (1+\mu^{-2}) \left\| \bsym{\varepsilon}_w + \dfrac{1}{2\sqrt{\mu}} \bsym{\varepsilon}_p \right\|_{0,\Omega}^2 \right) +  \dfrac{\beta_1^2 \mu}{4 K^2} \|\delta_u\|_{0,\Omega}^2 \right). \\
\end{split}
\label{eq:conv4.1.5}
\end{equation}
Therefore, we have
\begin{equation}
\| \delta_u \|_{0,\Omega}^2 
\leq \dfrac{C}{\mu}\left(\mu^{-1} \| \varepsilon_u \|_{0,\Omega}^2 + \mu^{-1} \| \bsym{\varepsilon}_p \|_{0,\Omega}^2 + (1+\mu^{-2}) \left\| \bsym{\varepsilon}_w + \dfrac{1}{2\sqrt{\mu}} \bsym{\varepsilon}_p \right\|_{0,\Omega}^2 \right).
\end{equation}
Finally, using the approximation properties \eqref{eq:approx_esdg}, we imply
\begin{equation}
\begin{split}
\| \delta_u \|_{0,\Omega}^2 \leq C(1+\mu^{-2}) h^{2(k+1)}.
\end{split}
\end{equation}
Using triangle inequality on $u - \widetilde{u}_h = \varepsilon_u + \delta_u$, we obtain our desired result.
\end{proof}
\end{theorem}
We remark that the error $u - \widetilde{u}_h$ consists of two parts. The difference $\delta_u$ between the numerical solution and the interpolation image depends on the viscosity coefficient $\mu$, while the interpolation error $\varepsilon_u$ does not. As we will see in our numerical results, the error does not vary significantly with the viscosity coefficient $\mu$.
}

\section{Numerical results}
\label{sec:num}
In this section, we illustrate some numerical examples. We carry out numerical experiments to see and compare the rates of convergence of the SDG method and the ESDG method. Polynomials with degree $k=1$ is used for SDG approximations. We are interested in the $L^2$ error of the unknown function $u$ and also that of the flux $\bsym{z} = \nabla u$. We recall the definitions in \eqref{eq:flux_esdg} and \eqref{eq:flux_sdg} for numerical approximations of the flux.

\begin{figure}[ht!]
\centering
\includegraphics[width=0.6\linewidth]{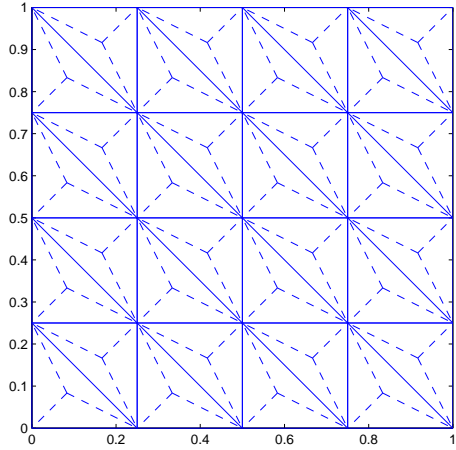}
\caption{An illustration of a staggered mesh on $\Omega = [0,1]^2$ with $N = 4$.}
\label{fig:mesh4}
\end{figure}

Throughout this section, we will take $\Omega = [0,1]^2 \subset \mathbb{R}^2$ and use a family of staggered meshes in all the experiments. We denote the number of uniform divisions in $[0,1]$ in the mesh by $N$. The domain $\Omega$ is partitioned into $N^2$ sub-squares with length $h = N^{-1}$. Each sub-square is then divided into two identical right-angled triangles by its diagonal. This constructs the initial triangulation $\mathcal{T}_u$. To construct the staggered mesh $\mathcal{T}$, we simply take the interior point $\nu$ as the centroid in each initial triangle $\mathcal{S}(\nu) \in \mathcal{T}_u$. Figure \ref{fig:mesh4} illustrates a member of this family with $N = 4$. Table \ref{tab:dof} compares the numbers of degrees of freedom in the discrete problems \eqref{eq:convdiff_sdg_system} for the SDG method and \eqref{eq:convdiff_esdg_system} for the ESDG method for this family of mesh. It can be seen that the number of degrees of freedom for the ESDG method is almost half of that of the SDG method for each level of mesh. Indeed, it can be worked out that for this family of mesh, the ratio of the numbers of the degrees of freedom is $7/12$ asymptotically:
\begin{equation}
\begin{split}
\text{dim}(U^h) & = 12N^2 + 4N,\\
\text{dim}(\widetilde{U}^h) & = 7N^2 + 2N + 1.
\end{split}
\end{equation}

\begin{table}[ht!]
\centering
\begin{tabular}{|c||c|c|}
\hline
$N$ & $\text{dim}(U^h)$ & $\text{dim}(\widetilde{U}^h)$ \\ 
\hline
2 & 56 & 33 \\
4 & 208 & 121 \\
8 & 800 & 465 \\
16 & 3136 & 1825 \\
32 & 12416 & 7233 \\
64 & 49408 & 28801 \\
\hline
\end{tabular}
\caption{Comparison of numbers of degrees of freedom: $\text{dim}(U^h)$ for SDG and $\text{dim} (\widetilde{U}^h)$ for ESDG.}
\label{tab:dof}
\end{table}

\subsection{Experiment 1: comparison to the EDG method}
\label{sec:exp1}

The purpose of this experiment is to compare our method with \cite{edg-cfd} in the same setting. In this experiment, the convection field $\bsym{b} = (b_1, b_2)$ is a constant vector.
The analytic solution of this experiment is given by
\begin{equation}
u(x,y) = xy\frac{(1-e^{b_1(x-1)})(1-e^{b_2(y-1)})}{(1-e^{b_1})(1-e^{b_2})}.
\label{eq:exp2_sol}
\end{equation}
For large values of $b_1$ and $b_2$, there is a boundary layer around the segments $x = 1$ and $y = 1$. 
The diffusivity $\mu$ are set to be $1$. The constant convection field is chosen to be $\bsym{b} = (20,20)$ and the problem is weakly convection-dominated. A homogeneous Dirichlet boundary condition $u \vert_{\partial \Omega} = 0$ is prescribed. The source function $f$ is computed accordingly. The SDG method \eqref{eq:convdiff_sdg} and the ESDG method \eqref{eq:convdiff_esdg} are used to solve the problem numerically. We examine the performance of the ESDG method by comparing the $L^2$ errors and the orders of $L^2$ convergence to the EDG method and the SDG method. 

\begin{figure}[ht!]
\centering
\includegraphics[width=0.6\linewidth]{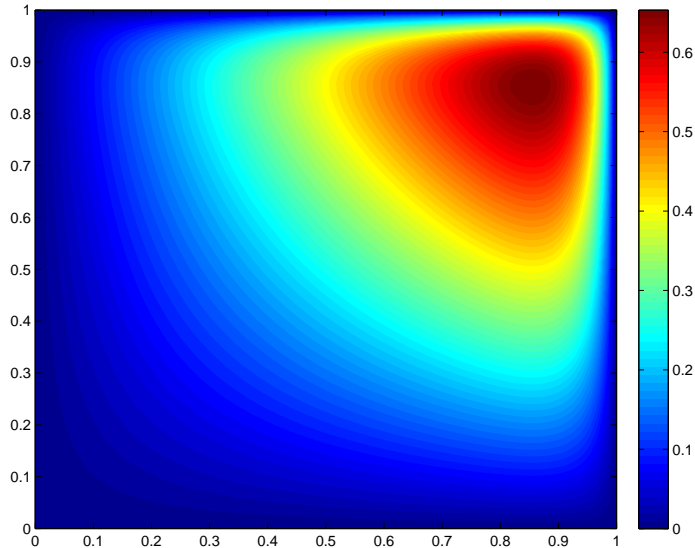}
\caption{A plot for the numerical solution $\widetilde{u}_h$ in Experiment 1.}
\label{fig:exp1plot}
\end{figure}

Figure \ref{fig:exp1plot} shows a plot of the numerical solution $\widetilde{u}_h$ of the ESDG method. Tables \ref{tab:1.Rn1} compare the convergence results of the SDG method and the ESDG method with various scales of diffusivity $\mu$. The second to the firth columns record the $L^2$ error and the orders of convergence of the potential and the flux for the SDG method. The sixth to the ninth columns record the $L^2$ error and the orders of convergence of the potential and the flux for the ESDG method. It can be seen that the approximated potential converge with an optimal order $2$ in $L^2$ error for both the SDG method and the ESDG method, which is the same for the EDG method in \cite{edg-cfd}. In particular, for $N = 32$ and $N = 64$, the ESDG method gives a $L^2$ error smaller that the SDG method and also the EDG method. However, similar to the EDG method for second-order elliptic problems \cite{edg-elliptic}, our numerical results show that the ESDG method only provides a suboptimal order $1$ of convergence of $L^2$ error of the approximated flux, while the SDG method provides an optimal order $2$. 

\begin{table}[ht!]
\centering
\begin{tabular}{|c||c c|c c||c c|c c|}
\hline
Mesh &  \multicolumn{2}{c|}{$\| u - u_h\|_{0, \Omega}$} &  \multicolumn{2}{c||}{$\| \bsym{z} - \bsym{z}_h\|_{0, \Omega}$} & \multicolumn{2}{c|}{$\| u - \widetilde{u}_h\|_{0, \Omega}$} &  \multicolumn{2}{c|}{$\| \bsym{z} - \widetilde{\bsym{z}}_h\|_{0, \Omega}$}\\ 
$N$ & Error & Order & Error & Order & Error & Order & Error & Order \\
\hline
    2 & 1.50e-01 &  -- & 2.48e+00 &  -- & 1.03e-01 &  -- & 2.31e+00 &  -- \\
    4 & 8.42e-02 & 0.84 & 1.63e+00 & 0.60 & 5.80e-02 & 0.82 & 2.00e+00 & 0.21 \\
    8 & 3.37e-02 & 1.32 & 6.96e-01 & 1.23 & 2.23e-02 & 1.38 & 1.38e+00 & 0.54 \\
   16 & 1.03e-02 & 1.72 & 2.15e-01 & 1.70 & 6.15e-03 & 1.86 & 7.92e-01 & 0.80 \\
   32 & 2.72e-03 & 1.91 & 5.71e-02 & 1.91 & 1.55e-03 & 1.99 & 4.13e-01 & 0.94 \\
   64 & 6.91e-04 & 1.98 & 1.45e-02 & 1.98 & 3.86e-04 & 2.00 & 2.09e-01 & 0.98 \\
\hline
\end{tabular}
\caption{History of convergence in Experiment 1.}
\label{tab:1.Rn1}
\end{table}

\subsection{Experiment 2: sensitivity of orders of convergence to diffusivity}
\label{sec:exp2}

\begin{figure}[ht!]
\centering
\includegraphics[width=0.6\linewidth]{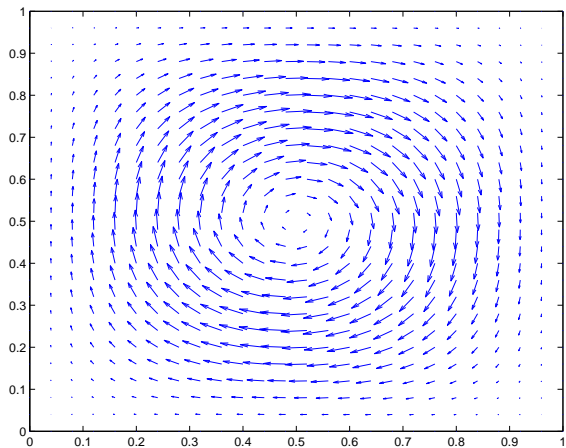}
\caption{The convection field $\bsym{b}$ in Experiment 2.}
\label{fig:exp2conv}
\end{figure}

The purpose of this experiment is to examine the performance of the ESDG method in terms of $L^2$ convergence and compare the ESDG method with the SDG method in various scales of diffusivity.
In this experiment, the convection field $\bsym{b} = (b_1, b_2)$ is set to be
\begin{equation}
\begin{split}
b_1 & = (1 - \cos(2 \pi x)) \sin(2 \pi y),\\
b_2 & = -\sin(2 \pi x) (1 - \cos(2 \pi y)).
\end{split}
\label{eq:exp2_conv}
\end{equation}
The analytic solution of this experiment is given by
\begin{equation}
u = \sin(2 \pi x) \cos(2 \pi y).
\label{eq:exp2_sol}
\end{equation}
Figure \ref{fig:exp2conv} shows a plot of the convection field $\bsym{b}$ in this experiment. We perform the experiment with different scales of diffusivity $\mu$. In particular, we are interested in observing the behaviour of the solutions when $\mu$ is small, i.e. the problem is convection-dominated. An inhomogeneous Dirichlet boundary condition is prescribed. The source function $f$ is computed accordingly. The SDG method \eqref{eq:convdiff_sdg} and the ESDG method \eqref{eq:convdiff_esdg} are used to solve the problem numerically. We examine the performance of the ESDG method by comparing the $L^2$ errors and the orders of $L^2$ convergence to the SDG method. 

Tables \ref{tab:2.Rn1}--\ref{tab:2.Rn10000} compare the convergence results of the SDG method and the ESDG method with various scales of diffusivity $\mu$. The second to the firth columns record the $L^2$ error and the orders of convergence of the potential and the flux for the SDG method. The sixth to the ninth columns record the $L^2$ error and the orders of convergence of the potential and the flux for the ESDG method. It can be seen that when the diffusivity $\mu$ is close to unity, the SDG method clearly outperforms the ESDG method. The convergence of the potential is optimal for both methods, while the convergence of the flux is optimal for the SDG method and suboptimal for the ESDG method. However, when the diffusivity $\mu$ reduces in scale, the optimal convergence of the flux for the SDG method is lost. By comparing the second column with the fourth column, and comparing the third column with the column, it can be seen that for convection-dominated situations, say $\mu \leq 10^{-3}$ in Tables \ref{tab:2.Rn1000}--\ref{tab:2.Rn10000}, the ESDG method has a comparable performance to the SDG method. With the considerable reduction in the size of the discrete problem, these results suggest that the ESDG method is favourable in convection-dominated situations. These observations in moderate problems and convection-dominated problems are in good agreement with the descriptions in \cite{edg-cfd}. \rev{Furthermore, we observe that, with a fixed mesh size, the $L^2$ error of the potential does not vary significantly with the viscosity coefficient $\mu$.}

\begin{table}[ht!]
\centering
\begin{tabular}{|c||c c|c c||c c|c c|}
\hline
Mesh &  \multicolumn{2}{c|}{$\| u - u_h\|_{0, \Omega}$} &  \multicolumn{2}{c||}{$\| \bsym{z} - \bsym{z}_h\|_{0, \Omega}$} & \multicolumn{2}{c|}{$\| u - \widetilde{u}_h\|_{0, \Omega}$} &  \multicolumn{2}{c|}{$\| \bsym{z} - \widetilde{\bsym{z}}_h\|_{0, \Omega}$}\\ 
$N$ & Error & Order & Error & Order & Error & Order & Error & Order \\
\hline
    2 & 3.28e-01 &  -- & 2.02e+00 &  -- & 6.37e-02 &  -- & 3.11e+00 &  -- \\
    4 & 1.23e-01 & 1.42 & 1.19e+00 & 0.76 & 1.05e-01 & -0.72 & 2.41e+00 & 0.37 \\ 
    8 & 3.58e-02 & 1.78 & 3.44e-01 & 1.79 & 3.81e-02 & 1.46 & 1.47e+00 & 0.71 \\
   16 & 9.31e-03 & 1.94 & 8.95e-02 & 1.94 & 1.06e-02 & 1.85 & 7.80e-01 & 0.92 \\
   32 & 2.35e-03 & 1.99 & 2.26e-02 & 1.99 & 2.72e-03 & 1.96 & 3.96e-01 & 0.98 \\
   64 & 5.89e-04 & 2.00 & 5.67e-03 & 2.00 & 6.85e-04 & 1.99 & 1.99e-01 & 0.99 \\
\hline
\end{tabular}
\caption{History of convergence for $\mu = 10^0$ in Experiment 2.}
\label{tab:2.Rn1}
\end{table}

\begin{table}[ht!]
\centering
\begin{tabular}{|c||c c|c c||c c|c c|}
\hline
Mesh &  \multicolumn{2}{c|}{$\| u - u_h\|_{0, \Omega}$} &  \multicolumn{2}{c||}{$\| \bsym{z} - \bsym{z}_h\|_{0, \Omega}$} & \multicolumn{2}{c|}{$\| u - \widetilde{u}_h\|_{0, \Omega}$} &  \multicolumn{2}{c|}{$\| \bsym{z} - \widetilde{\bsym{z}}_h\|_{0, \Omega}$}\\ 
$N$ & Error & Order & Error & Order & Error & Order & Error & Order \\
\hline
    2 & 1.07e+00 &  -- & 1.19e+01 &  -- & 6.40e-01 &  -- & 9.39e+00 &  -- \\
    4 & 2.16e-01 & 2.31 & 3.89e+00 & 1.61 & 2.62e-01 & 1.29 & 6.11e+00 & 0.62 \\
    8 & 5.58e-02 & 1.95 & 1.49e+00 & 1.38 & 5.55e-02 & 2.24 & 2.60e+00 & 1.23 \\
   16 & 1.39e-02 & 2.01 & 4.96e-01 & 1.59 & 1.24e-02 & 2.16 & 1.02e+00 & 1.35 \\
   32 & 3.38e-03 & 2.03 & 1.43e-01 & 1.79 & 2.88e-03 & 2.11 & 4.35e-01 & 1.23 \\
   64 & 8.33e-04 & 2.02 & 3.79e-02 & 1.92 & 7.01e-04 & 2.04 & 2.04e-01 & 1.09 \\
\hline
\end{tabular}
\caption{History of convergence for $\mu = 10^{-2}$ in Experiment 2.}
\label{tab:2.Rn100}
\end{table}

\begin{table}[ht!]
\centering
\begin{tabular}{|c||c c|c c||c c|c c|}
\hline
Mesh &  \multicolumn{2}{c|}{$\| u - u_h\|_{0, \Omega}$} &  \multicolumn{2}{c||}{$\| \bsym{z} - \bsym{z}_h\|_{0, \Omega}$} & \multicolumn{2}{c|}{$\| u - \widetilde{u}_h\|_{0, \Omega}$} &  \multicolumn{2}{c|}{$\| \bsym{z} - \widetilde{\bsym{z}}_h\|_{0, \Omega}$}\\ 
$N$ & Error & Order & Error & Order & Error & Order & Error & Order \\
\hline
    2 & 1.91e+00 &  -- & 3.03e+01 &  -- & 1.35e+00 &  -- & 2.16e+01 &  -- \\
    4 & 3.11e-01 & 2.62 & 8.52e+00 & 1.83 & 4.37e-01 & 1.63 & 1.36e+01 & 0.66 \\
    8 & 8.17e-02 & 1.93 & 3.93e+00 & 1.12 & 8.27e-02 & 2.40 & 4.91e+00 & 1.47 \\
   16 & 2.07e-02 & 1.98 & 1.57e+00 & 1.32 & 1.82e-02 & 2.18 & 1.74e+00 & 1.50 \\
   32 & 4.78e-03 & 2.11 & 5.35e-01 & 1.55 & 4.15e-03 & 2.13 & 6.71e-01 & 1.37 \\
   64 & 1.09e-03 & 2.13 & 1.58e-01 & 1.76 & 9.71e-04 & 2.10 & 2.67e-01 & 1.33 \\
\hline
\end{tabular}
\caption{History of convergence for $\mu = 2 \times 10^{-3}$ in Experiment 2.}
\label{tab:2.Rn500}
\end{table}

\begin{table}[ht!]
\centering
\begin{tabular}{|c||c c|c c||c c|c c|}
\hline
Mesh &  \multicolumn{2}{c|}{$\| u - u_h\|_{0, \Omega}$} &  \multicolumn{2}{c||}{$\| \bsym{z} - \bsym{z}_h\|_{0, \Omega}$} & \multicolumn{2}{c|}{$\| u - \widetilde{u}_h\|_{0, \Omega}$} &  \multicolumn{2}{c|}{$\| \bsym{z} - \widetilde{\bsym{z}}_h\|_{0, \Omega}$}\\ 
$N$ & Error & Order & Error & Order & Error & Order & Error & Order \\
\hline
    2 & 2.33e+00 &  -- & 4.23e+01 &  -- & 2.42e+00 &  -- & 3.98e+01 &  -- \\
    4 & 3.61e-01 & 2.69 & 1.10e+01 & 1.95 & 6.06e-01 & 2.00 & 2.09e+01 & 0.93 \\
    8 & 9.92e-02 & 1.86 & 5.66e+00 & 0.95 & 1.09e-01 & 2.47 & 7.20e+00 & 1.54 \\
   16 & 2.58e-02 & 1.94 & 2.37e+00 & 1.26 & 2.30e-02 & 2.25 & 2.41e+00 & 1.58 \\
   32 & 5.99e-03 & 2.11 & 8.87e-01 & 1.42 & 4.92e-03 & 2.23 & 8.58e-01 & 1.49 \\
   64 & 1.31e-03 & 2.20 & 2.77e-01 & 1.68 & 1.14e-03 & 2.11 & 3.34e-01 & 1.36 \\
\hline
\end{tabular}
\caption{History of convergence for $\mu = 10^{-3}$ in Experiment 2.}
\label{tab:2.Rn1000}
\end{table}

\begin{table}[ht!]
\centering
\begin{tabular}{|c||c c|c c||c c|c c|}
\hline
Mesh &  \multicolumn{2}{c|}{$\| u - u_h\|_{0, \Omega}$} &  \multicolumn{2}{c||}{$\| \bsym{z} - \bsym{z}_h\|_{0, \Omega}$} & \multicolumn{2}{c|}{$\| u - \widetilde{u}_h\|_{0, \Omega}$} &  \multicolumn{2}{c|}{$\| \bsym{z} - \widetilde{\bsym{z}}_h\|_{0, \Omega}$}\\ 
$N$ & Error & Order & Error & Order & Error & Order & Error & Order \\
\hline
    2 & 2.78e+00 &  -- & 5.53e+01 &  -- & 4.65e+00 &  -- & 7.80e+01 &  -- \\
    4 & 4.34e-01 & 2.68 & 1.39e+01 & 2.00 & 9.42e-01 & 2.31 & 3.34e+01 & 1.22 \\
    8 & 1.22e-01 & 1.83 & 7.92e+00 & 0.81 & 1.67e-01 & 2.50 & 1.15e+01 & 1.53 \\
   16 & 3.20e-02 & 1.93 & 3.37e+00 & 1.23 & 3.13e-02 & 2.41 & 3.77e+00 & 1.61 \\
   32 & 7.69e-03 & 2.06 & 1.39e+00 & 1.28 & 6.15e-03 & 2.35 & 1.20e+00 & 1.65 \\
   64 & 1.66e-03 & 2.21 & 4.72e-01 & 1.56 & 1.33e-03 & 2.20 & 4.27e-01 & 1.49 \\
\hline
\end{tabular}
\caption{History of convergence for $\mu = 5 \times 10^{-4}$ in Experiment 2.}
\label{tab:2.Rn2000}
\end{table}

\begin{table}[ht!]
\centering
\begin{tabular}{|c||c c|c c||c c|c c|}
\hline
Mesh &  \multicolumn{2}{c|}{$\| u - u_h\|_{0, \Omega}$} &  \multicolumn{2}{c||}{$\| \bsym{z} - \bsym{z}_h\|_{0, \Omega}$} & \multicolumn{2}{c|}{$\| u - \widetilde{u}_h\|_{0, \Omega}$} &  \multicolumn{2}{c|}{$\| \bsym{z} - \widetilde{\bsym{z}}_h\|_{0, \Omega}$}\\ 
$N$ & Error & Order & Error & Order & Error & Order & Error & Order \\
\hline
    2 & 3.78e+00 &  -- & 8.05e+01 &  -- & 1.14e+01 &  -- & 1.94e+02 &  -- \\
    4 & 6.29e-01 & 2.59 & 2.08e+01 & 1.96 & 2.05e+00 & 2.48 & 7.21e+01 & 1.43 \\
    8 & 1.66e-01 & 1.93 & 1.20e+01 & 0.79 & 3.41e-01 & 2.59 & 2.40e+01 & 1.59 \\
   16 & 4.28e-02 & 1.95 & 5.03e+00 & 1.26 & 5.65e-02 & 2.59 & 7.94e+00 & 1.60 \\
   32 & 1.07e-02 & 2.00 & 2.26e+00 & 1.15 & 9.58e-03 & 2.56 & 2.28e+00 & 1.80 \\
   64 & 2.39e-03 & 2.16 & 8.86e-01 & 1.35 & 1.79e-03 & 2.42 & 6.88e-01 & 1.73 \\
\hline
\end{tabular}
\caption{History of convergence for $\mu = 2 \times 10^{-4}$ in Experiment 2.}
\label{tab:2.Rn5000}
\end{table}

\begin{table}[ht!]
\centering
\begin{tabular}{|c||c c|c c||c c|c c|}
\hline
Mesh &  \multicolumn{2}{c|}{$\| u - u_h\|_{0, \Omega}$} &  \multicolumn{2}{c||}{$\| \bsym{z} - \bsym{z}_h\|_{0, \Omega}$} & \multicolumn{2}{c|}{$\| u - \widetilde{u}_h\|_{0, \Omega}$} &  \multicolumn{2}{c|}{$\| \bsym{z} - \widetilde{\bsym{z}}_h\|_{0, \Omega}$}\\ 
$N$ & Error & Order & Error & Order & Error & Order & Error & Order \\
\hline
    2 & 5.33e+00 &  -- & 1.20e+02 &  -- & 2.28e+01 &  -- & 3.87e+02 &  -- \\
    4 & 8.87e-01 & 2.59 & 3.06e+01 & 1.98 & 3.91e+00 & 2.54 & 1.38e+02 & 1.49 \\
    8 & 2.21e-01 & 2.01 & 1.64e+01 & 0.90 & 6.20e-01 & 2.66 & 4.37e+01 & 1.66 \\
   16 & 5.43e-02 & 2.02 & 6.75e+00 & 1.28 & 1.03e-01 & 2.60 & 1.50e+01 & 1.54 \\
   32 & 1.39e-02 & 1.97 & 3.10e+00 & 1.12 & 1.52e-02 & 2.76 & 4.10e+00 & 1.87 \\
   64 & 3.17e-03 & 2.13 & 1.32e+00 & 1.23 & 2.52e-03 & 2.59 & 1.15e+00 & 1.84 \\
\hline
\end{tabular}
\caption{History of convergence for $\mu = 10^{-4}$ in Experiment 2.}
\label{tab:2.Rn10000}
\end{table}

\subsection{Experiment 3: uniform stability in $L^2$ energy with respect to diffusivity}
\label{sec:exp3}

The purpose of this experiment is to examine the stability in $L^2$ energy of the ESDG method in various scales of diffusivity. We first observe that it is actually possible to derive different discretizations for the convection term and the diffusion term, and we will compare our skew-symmetric discretization with two types of non-skew-symmetric discretizations. Given $\theta \in [0,1]$. we modify the definitions of the auxiliary variables in \eqref{eq:newvar} by
\begin{equation}
\begin{split}
\bsym{w} & = \sqrt{\mu} \nabla u - \dfrac{\theta}{\sqrt{\mu}} \bsym{b}u, \\
\bsym{z} & = \bsym{b}u.
\end{split}
\end{equation}
Then we can use the same idea as \eqref{eq:convdiff_esdg} to obtain a new method. The discrete convection term is then modified accordingly as
\begin{equation}
\bsym{b} \cdot \nabla_h = -\theta \widetilde{B} M^{-1} \widetilde{R}^T + (1 - \theta) \widetilde{R} M^{-1} \widetilde{B}^T.
\end{equation}
In particular, when $\theta = 1/2$, it is reduced to ESDG method \eqref{eq:convdiff_esdg} with a skew-symmetric discretization of the convection term proposed in Section \ref{sec:esdg}. We will compare the discretizations with $\theta = 0$, $\theta = 1/2$ and $ \theta = 1$, and observe the advantages brought by the spectro-consistent discretization with the novel splitting of the convection term and the diffusion term. We remark that a similar experiment is performed on the SDG method for incompressible Navier-Stokes equations in \cite{sdg-ns1}. 

In this experiment, the convection field $\bsym{b} = (b_1, b_2)$ is identical to Experiment 2.
\begin{equation}
\begin{split}
b_1 & = (1 - \cos(2 \pi x)) \sin(2 \pi y),\\
b_2 & = -\sin(2 \pi x) (1 - \cos(2 \pi y)). \\
\end{split}
\label{eq:exp3_conv}
\end{equation}
The analytic solution of this experiment is given by
\begin{equation}
u = \sin(2 \pi x) \sin(2 \pi y).
\label{eq:exp3_sol}
\end{equation}
We perform the experiment with different scales of diffusivity $\mu$. In particular, we are interested in observing the behaviour of the solutions when $\mu$ is small, i.e. the problem is convection-dominated. A homogeneous Dirichlet boundary condition $u \vert_{\partial \Omega} = 0$ is prescribed. The source function $f$ is computed accordingly. The ESDG method \eqref{eq:convdiff_esdg} is used to solve the problem numerically. We use a mesh with size $N = 32$. We are interested in the $L^2$ norm $\| \widetilde{\bsym{z}}_h \|_{0,\Omega}$ of the approximation $\widetilde{\bsym{z}}_h$ of the flux $\bsym{z}$. By a direct computation, it is easy to see that $\| \bsym{z} \|_{0,\Omega} = \sqrt{2}\pi \approx 4.4429$.

Tables \ref{tab:3.L2z} records the $L^2$ norm $\| \widetilde{\bsym{z}}_h \|_{0,\Omega}$ of the approximation $\widetilde{\bsym{z}}_h$ with the three different discretizations. For more moderate problems $\mu > 10^{-3}$, it can be seen that all the three discretizations provide a approximation $\widetilde{\bsym{z}}_h$ with the $L^2$ norm close to the value $\sqrt{2}\pi$. However, for convection dominated problems, the skew symmetric discretization $\theta = 1/2$ clearly outperforms the other two discretizations. In spite of the machine error due to an ill-conditioned linear system as the diffusivity tends to zero, the $L^2$ norm $\| \widetilde{\bsym{z}}_h \|_{0,\Omega}$ of the approximation $\widetilde{\bsym{z}}_h$ is around a constant when $\theta = 1/2$. Meanwhile, for the other two discretizations, the $L^2$ norm $\| \widetilde{\bsym{z}}_h \|_{0,\Omega}$ of the approximation $\widetilde{\bsym{z}}_h$ blows up as the diffusivity tends to zero.

\begin{table}[ht!]
\centering
\begin{tabular}{|r||c|c|c|}
\hline
Diffusivity &  \multicolumn{3}{c|}{$\| z_h\|_{0, \Omega}$}\\ \cline{2-4}
$\mu$ & $\theta = 0$ & $\theta = 1/2$ & $\theta = 1$\\
\hline
$10^{0}$ & 4.43e+00 & 4.43e+00 & 4.43e+00 \\
$10^{-2}$ & 4.47e+00 & 4.47e+00 & 4.47e+00 \\
$2 \times 10^{-3}$ & 4.48e+00 & 4.49e+00 & 4.59e+00 \\
$10^{-3}$ & 6.33e+00 & 4.52e+00 & 9.31e+00 \\
$5 \times 10^{-4}$ & 1.12e+03 & 4.59e+00 & 2.01e+03 \\
$2 \times 10^{-4}$ & 1.81e+03 & 4.88e+00 & 8.73e+02 \\
$10^{-4}$ & 1.55e+05 & 5.52e+00 & 8.51e+04 \\
\hline
\end{tabular}
\caption{Record of $\| \bsym{z}_h \|_{0,\Omega}$ in Experiment 3.}
\label{tab:3.L2z}
\end{table}

\section{Conclusion}
\label{sec:conclusion}
In this paper, we develop an embedded staggered discontinuous Galerkin method for the convection-diffusion equation. Thanks to the design of the SDG finite element spaces, the new method provides local and global conservations, and does not require the introduction of carefully designed stabilization terms or flux conditions. Furthermore, $L^2$ stability is achieved by a skew-symmetric discretization of the convection term. Numerical results are presented to show the robustness of the method with respect to diffusivity. On the other hand, the method seeks reduced approximations in a subspace of the SDG finite element space. In convection-dominated problems, like other DG methods, the convergence are optimal in potential and suboptimal in flux, as our numerical results have shown.


\begin{thebibliography}{99}

\bibitem{unified-dg}
{\sc D. N. Arnold, F. Brezzi, B. Cockburn, and L. D. Marini}, 
\emph{Unified analysis of discontinuous Galerkin methods for elliptic problems}, 
SIAM J. Numer. Anal., 39 (2002), pp. 1749--1779.

\bibitem{sfvm}
{\sc P. Blanc, R. Eymard, R. Herbin},
\emph{A staggered finite volume scheme on general meshes for the generalized Stokes problem in two space dimensions},
International Journal on Finite Volumes, 2 (2005), pp. 1--31.

\bibitem{sdm}
{\sc B. J. Boersma},
\emph{A staggered compact finite difference formulation for the compressible Navier-Stokes equations},
J. Comput. Phys., 208 (2005), pp. 675--690.

\bibitem{braess}
{\sc D. Braess},
\emph{Finite elements. Theory, fast solvers, and applications in elasticity theory},
Cambridge University Press, Cambridge, 2007.

\bibitem{buffa06}
{\sc A. Buffa, T. J. R. Hughes, G. Sangalli}, 
\emph{Analysis of a multiscale discontinuous Galerkin method for convection-diffusion problems}, 
SIAM J. Numer. Anal., 44 (2006), pp. 1420--1440.

\bibitem{carrero05}
{\sc J. Carrero, B. Cockburn, D. Sch\"{o}tzau},
\emph{Hybridized globally divergence-free LDG methods. Part I: The Stokes problem},
Math. Comput., 75 (2005), pp. 533--563.

\bibitem{sdg-ns1}
{\sc S. W. Cheung, E. Chung, H. H. Kim, Y. Qian},
\emph{Staggered discontinuous Galerkin methods for incompressible Navier-Stokes equations},
J. Comput. Phys., 302 (2015), pp. 251--266.

\bibitem{sdg-ibm}
{\sc S. W. Cheung, E. Chung, H. H. Kim},
\emph{Staggered discontinuous Galerkin approximation for immersed boundary method},
arXiv preprint, arXiv:1609.01046, 2016.

\bibitem{semi}
{\sc E. T. Chung, Q. Du, J. Zou},
\emph{Convergence analysis on a finite volume method for Maxwell's equations in non-homogeneous media},
SIAM J. Numer. Anal., 41 (2003), pp. 37--63.

\bibitem{fully}
{\sc E. T. Chung, B. Engquist},
\emph{Convergence analysis of fully discrete finite volume methods for Maxwell's equations in nonhomogeneous media},
SIAM J. Numer. Anal., 43 (2005), pp. 303--317.

\bibitem{newdg}
{\sc E. T. Chung, B. Engquist},
\emph{Optimal discontinuous Galerkin methods for wave propagation},
SIAM J. Numer. Anal., 44 (2006), pp. 2131--2158.

\bibitem{newdg1}
{\sc E. T. Chung, B. Engquist},
\emph{Optimal discontinuous Galerkin methods for the acoustic wave equation in higher dimensions},
SIAM J. Numer. Anal., 47 (2009), pp. 3820--3848.

\bibitem{meta}
{\sc E. T. Chung, P. Ciarlet},
\emph{A staggered discontinuous Galerkin method for wave propagation in media with dielectrics and meta-materials},
J. Comput. Appl. Math., 239 (2013), pp. 189--207.

\bibitem{jcp-max}
{\sc E. T. Chung, P. Ciarlet, T. F. Yu},
\emph{Convergence and superconvergence of staggered discontinuous Galerkin methods for the three-dimensional Maxwell's equations on Cartesian grids},
J. Comput. Phys., 235 (2013), pp. 14--31.

\bibitem{sdg-hdg}
{\sc E. Chung, B. Cockburn, G. Fu},
\emph{The staggered DG method is the limit of a hybridizable DG method},
SIAM J. Numer. Anal., 52 (2014), pp. 915--932.


\bibitem{sdg-hdg1}
{\sc E. Chung, B. Cockburn, G. Fu},
\emph{ The staggered DG method is the limit of a hybridizable DG method. Part II: The Stokes flow.},
J. Sci. Comput., 66 (2016), pp. 870-887.

\bibitem{geo}
{\sc E. T. Chung, C. Y. Lam and J. Qian},
\emph{A staggered discontinuous Galerkin method for the simulation of seismic waves with surface topography},
Geophysics, 80 (2015), pp. T119-T135.




\bibitem{curlcurl}
{\sc E. T. Chung, C. S. Lee},
\emph{A staggered discontinuous Galerkin method for the curl-curl operator},
 IMA J. Numer. Anal., 32 (2012), pp. 1241--1265.

\bibitem{convdiff}
{\sc E. T. Chung, C. S. Lee},
\emph{A staggered discontinuous Galerkin method for the convection-diffusion equation},
J. Numer. Math., 20 (2012), pp. 1--31.

\bibitem{sdg-ns2}
{\sc E. T. Chung, W. Qiu},
\emph{Analysis of a SDG method for the incompressible Navier-Stokes equations},
Submitted.

\bibitem{ciarlet}
{\sc P. Ciarlet},
\emph{The Finite Element Method for Elliptic
Problems},
 North-Holland, Amsterdam.

\bibitem{hdg}
{\sc B. Cockburn, J. Gopalakrishnan, N.C. Nguyen, J. Peraire, F.-J. Sayas},
\emph{Analysis of an HDG method for Stokes flow},
Math. Comput., 80 (2011), pp. 723--760.

\bibitem{unified-hdg}
{\sc B. Cockburn, J. Gopalakrishnan, R. Lazarov}, 
\emph{Unified hybridization of discontinuous Galerkin, mixed and continuous Galerkin methods for second order elliptic problems}, 
SIAM J. Numer. Anal. 47 (2009) pp. 1319--1365.


\bibitem{cockburn05}
{\sc B. Cockburn, G. Kanschat, D. Sch\"{o}tzau},
\emph{A locally conservative LDG method for the incompressible Navier-Stokes equations},
Math. Comp., 74 (2005), pp. 1067--1095.

\bibitem{hdg-convdiff}
{\sc B. Cockburn, B. Dong, J. Guzman, M. Restelli, R. Sacco}, 
\emph{A hybridizable discontinuous Galerkin method for steady-state convection-diffusion-reaction
problems}, 
SIAM J. Sci. Comput., 31 (2009), 3827--3846.

\bibitem{edg-elliptic}
{\sc B. Cockburn, J. Guzmán, S.-C. Soon, H.K. Stolarski}, 
\emph{An analysis of the embedded discontinuous Galerkin method for second-order elliptic problems},
SIAM J. Numer. Anal., 47 (2009) pp. 2686--2707.

\bibitem{ldg}
{\sc B. Cockburn, G. Kanschat, D. Sch\"{o}tzau, C. Schwab},
\emph{Local discontinuous Galerkin methods for the Stokes system},
SIAM J. Numer. Anal., 40 (2002), pp. 319--343.

\bibitem{conv-diff}
{\sc B. Cockburn, C.-W. Shu},
\emph{The local discontinuous Galerkin method for time-dependent convection-diffusion systems},
SIAM J. Numer. Anal., 35 (1998), 2440--2463.

\bibitem{edg-les}
{\sc P. Fernandez, N. C. Nguyen, X. Roca, J. Peraire},
\emph{Implicit large-eddy simulation of compressible flows using the Interior Embedded Discontinuous Galerkin method},
54th AIAA Aerospace Sciences Meeting, AIAA SciTech, (AIAA 2016-1332)

\bibitem{edg-shell}
{\sc S. G\"{u}zey, B. Cockburn, and H.K. Stolarski}, 
\emph{The embedded discontinuous Galerkin methods: Application to linear shells problems}, 
Internat. J. Numer. Methods Engrg., 70 (2007), 757--790.

\bibitem{sdm1}
{\sc F. H. Harlow, J. E. Welch},
\emph{Numerical calculation of time-dependent viscous incompressible flow of fluid with a free surface},
Phys. Fluids, 8 (1965), pp. 2182--2189.

\bibitem{houston09}
{\sc P. Houston, D. Sch\"{o}tzau, X. Wei},
\emph{A mixed DG method for linearized incompressible magnetohydrodynamics},
  J. Sci. Comp., 40 (2009), pp. 281--314.

\bibitem{hughes06}
{\sc J. T. R. Hughes, G. Scovazzi, P. B. Bochev, A. Buffa}, 
\emph{A multiscale discontinuous Galerkin method with the computational structure of a continuous Galerkin method}, 
Comput. Methods Appl. Mech. Engrg., 195 (2006), pp. 2761--2787.

\bibitem{kcl-2013-Stokes-SDG}
{\sc H. H. Kim, E. T. Chung, C. S. Lee},
\emph{A staggered discontinuous Galerkin method for the Stokes system},
SIAM J. Numer. Anal., 51 (2013), pp. 3327--3350.

\bibitem{kcl-2013-FETI-DP-Stokes}
{\sc H. H. Kim, E. T. Chung, C. S. Lee},
\emph{{FETI}-{DP} preconditioners for a staggered discontinuous {G}alerkin formulation
of the two-dimensional Stokes problem},
Comput. \& Math. Appl., 68 (2014), pp. 2233-2250.

\bibitem{Liu-Shu}
{\sc J.-G. Liu, C.-W. Shu},
\emph{A high-order discontinuous Galerkin method for 2D incompressible flows},
J. Comput. Phys., 160 (2000), pp. 577--596.

\bibitem{nguyen11}
{\sc N. C. Nguyen, J. Peraire, B. Cockburn},
\emph{An implicit high-order hybridizable discontinuous Galerkin method for the incompressible Navier-Stokes equations},
J. Comput. Phys., 230 (2011), pp. 1147--1170.

\bibitem{edg-cfd}
{\sc N.C. Nguyen, J. Peraire, B. Cockburn}, 
\emph{A class of embedded discontinuous Galerkin methods for computational fluid dynamics},
J. Comput. Phys., 302 (2015), pp. 674--692.

\bibitem{edg-ns}
{\sc J. Peraire, N.C. Nguyen, B. Cockburn}, 
\emph{An embedded discontinuous Galerkin method for the compressible Euler and Navier–Stokes equations },
(AIAA Paper 2011-3228), in: Proceedings of the 20th AIAA Computational Fluid Dynamics Conference, Honolulu, HI, 2011.

\bibitem{first-dg}
{\sc W. H. Reed, T. R. Hill}, 
\emph{Triangular Mesh Methods for the Neutron Transport Equation}, 
Tech. Report LA-UR-73-479, Los Alamos Scientific Laboratory, Los Alamos, NM, 1973.

\bibitem{hpdg}
{\sc D. Sch\"{o}tzau, C. Schwab, A. Toselli},
\emph{Mixed hp-DGFEM for incompressible flows},
  SIAM J. Numer. Anal., 40 (2003), pp. 2171--2194.

\bibitem{shahbazi07}
{\sc K. Shahbazi, P. F Fischer, C.R. Ethier},
\emph{A high-order discontinuous Galerkin method for the unsteady incompressible Navier-Stokes equations},
J. Comput. Phys. 222 (2007), pp. 391--407.

\end{thebibliography}
\end{document}